\documentclass[11pt]{article}

\usepackage{titlesec}
\usepackage{a4wide}
\usepackage{graphicx}
\usepackage{float}
\usepackage{amssymb}
\usepackage{amsmath}
\usepackage{amsthm}
\usepackage{appendix}
\usepackage{color}
\usepackage{array}
\usepackage{pstricks}
\usepackage{eucal}
\usepackage{mathrsfs}
\usepackage{esint}
\usepackage{stmaryrd}

\newtheorem{theorem}{Theorem}[section]
\newtheorem{lemma}{Lemma}[section]
\newtheorem{remark}{Remark}[section]

\newtheorem{corollary}{Corollary}[section]
\newtheorem{proposition}{Proposition}[section]

\newcommand{\nc}{\newcommand}
\nc{\dsp}{\displaystyle}
\nc{\path}[1]{../Figures/#1}
\nc{\llbr}{\llbracket}
\nc{\rrbr}{\rrbracket}
\nc{\lbr}{\lbrack}
\nc{\rbr}{\rbrack}
\nc{\eps}{\delta}
\nc{\mrm}[1]{\mathrm{#1}}
\nc{\Tau}{\mathbb{T}}
\nc{\vp}{\varphi}
\nc{\axi}[1]{\breve{#1}}
\nc{\bfx}{\boldsymbol{x}}
\nc{\bfy}{\boldsymbol{y}}
\nc{\bfn}{\mathbf{n}}

\nc{\R}{\mathbb{R}}
\nc{\C}{\mathbb{C}}
\nc{\T}{\mrm{T}}
\nc{\mL}{\mrm{L}}
\nc{\mH}{\mrm{H}}
\nc{\mX}{\mrm{X}}
\nc{\bfH}{\boldsymbol{H}}
\nc{\Imag}{\Im \mrm{m}}
\nc{\Rg}{\mathcal{R}}
\nc{\Curl}{\mathbf{curl}}
\nc{\curl}{\mathrm{curl}}
\nc{\bfL}{\boldsymbol{L}}
\nc{\bsb}[1]{\boldsymbol{#1}}

\nc{\bfU}{\mathbf{U}}
\nc{\bfV}{\mathbf{V}}
\nc{\bfW}{\mathbf{W}}
\nc{\bfF}{\mathbf{F}}
\nc{\bfu}{\mathbf{u}}
\nc{\bfv}{\mathbf{v}}
\nc{\bfp}{\mathbf{p}}
\nc{\bfq}{\mathbf{q}}
\nc{\bfw}{\mathbf{w}}
\nc{\bft}{\boldsymbol{t}}
\nc{\bfh}{\boldsymbol{h}}
\nc{\bfe}{\boldsymbol{e}}
\nc{\bfalpha}{\boldsymbol{\alpha}}
\nc{\bfbeta}{\boldsymbol{\beta}}

\nc{\dir}{\textsc{d}}
\nc{\neu}{\textsc{n}}
\nc{\metal}{\tiny\Sigma}
\nc{\inc}{\mrm{inc}}
\nc{\loc}{\mrm{loc}}
\nc{\comp}{\mrm{comp}}
\nc{\bcdir}{\Lambda_{\dir}}
\nc{\bcneu}{\Lambda_{\neu}}

\renewcommand{\div}{\mrm{div}}
\nc{\vect}[1]{\boldsymbol{#1}}

\nc{\vrh}{\varrho}
\nc{\vth}{\vartheta}
\nc{\lbd}{\lambda}
\nc{\mbX}{\mathbb{X}}
\nc{\mbY}{\mathbb{Y}}
\nc{\mbH}{\mathbb{H}}
\nc{\mbL}{\mathbb{L}}

\nc{\SL}{\boldsymbol{\mrm{SL}}}
\nc{\DL}{\boldsymbol{\mrm{DL}}}
\nc{\Green}{\mrm{G}}
\nc{\Cauchy}{\mathcal{C}}
\nc{\Lift}{\boldsymbol{\Phi}}

\nc{\step}{\mrm{h}}
\nc{\wh}{\widehat}

\nc{\TrSp}{\boldsymbol{H}^{-1/2}_{\times,\div}}
\nc{\SgTrSp}{\mathbb{X}}
\nc{\SplitZ}{\boldsymbol{\mathcal{Z}}}
\nc{\SplitN}{\boldsymbol{\mathcal{N}}}

\nc{\FigUn}{\textbf{Fig.1}\,}
\nc{\FigDe}{\textbf{Fig.2}\,}
\nc{\FigTr}{\textbf{Fig.3}\,}

\nc{\TrD}{\boldsymbol{\gamma}_{\dir}}
\nc{\TrT}{\boldsymbol{\gamma}_{\bft}}
\nc{\TrN}{\boldsymbol{\gamma}_{\neu}}
\nc{\TrDc}{\boldsymbol{\gamma}_{\dir,c}}
\nc{\TrNc}{\boldsymbol{\gamma}_{\neu,c}}
\nc{\trd}{\gamma_{\mrm{d}}}
\nc{\Tr}{\boldsymbol{\gamma}}
\nc{\Trc}{\boldsymbol{\gamma}_{c}}

\nc{\Normal}{\mathbf{n}}
\nc{\pll}{\parallel}
\nc{\SM}{\mrm{SM}}
\nc{\mrB}{\mrm{B}}
\nc{\mrS}{\mrm{S}}
\nc{\mrD}{\mrm{D}}
\nc{\mrT}{\mrm{T}}

\nc{\TrOp}{\mathcal{T}}
\nc{\mrA}{\mrm{A}}
\nc{\mrK}{\mrm{K}}
\nc{\mrQ}{\mrm{Q}}
\nc{\mrR}{\mrm{R}}
\nc{\mrId}{\mrm{Id}}
\nc{\DTN}{\mrm{T}}

\nc{\Range}{\mathcal{R}}
\nc{\Ker}{\mrm{ker}}

\nc{\llb}{\llbracket}
\nc{\rrb}{\rrbracket}
\nc{\ctr}{\mathfrak}
\nc{\ctru}{\mathfrak{u}}
\nc{\ctrv}{\mathfrak{v}}

\nc{\kernel}{\mathscr{G}}
\nc{\mrint}{\mathrm{in}}
\nc{\mrout}{\mathrm{out}}

\definecolor{green2}{rgb}{0,0.5,0}

\date{}
\title{\textbf{Essential spectrum  of  local\\  multi-trace boundary integral operators}}
\author{
X.~Claeys\footnote{Laboratoire Jacques-Louis Lions, UPMC University of Paris 6 and CNRS UMR 7598, 75005, Paris, France}
$\!\,^{,}$\footnote{INRIA-Paris-Rocquencourt, EPC Alpines, Le Chesnay Cedex, France}
}

\begin{document}

\maketitle
\begin{abstract} 
  Considering pure transmission scattering problems in piecewise constant media,
  we derive an exact analytic formula for the spectrum of the corresponding
  local multi-trace boundary integral operators in the case where the geometrical
  configuration does not involve any junction point and all wave numbers equal.
  We deduce from this the essential spectrum in the case where wave numbers vary.
  Numerical evidences of these theoretical results are also presented.
\end{abstract}

\section*{Introduction}

Many applications involve the simulation of wave propagation phenomena in media with 
piecewise constant material characteristics that can be modeled by elliptic 
partial differential equations with piecewise constant coefficients. In such 
situations, the computational domain is naturally partitioned into sub-domains 
corresponding to the values of material characteristics. 

As regards numerical approaches to be used to tackle wave propagation problems, 
although one may opt for finite elements or similar volume methods, boundary integral 
equation methods provide accurate alternatives that are less prone to such 
undesirable effects as numerical dispersion. Admittedly discretization of 
boundary integral equations leads to dense ill conditioned matrices which raises
numerical challenges and requires careful implementation, but many progresses achieved 
in the past decade  (efficient preconditioners, fast multipole methods, adapted 
quadrature techniques) now place boundary integral equation techniques as a serious 
numerical alternative for high performance computations.

In the context of parallel computing, it becomes desirable to embed integral 
equation methods into a domain decomposition paradigm. The Boundary Element Tearing 
and Interconnecting method (BETI) was developed in this spirit, more than a decade ago, as an 
integral equation counterpart of the FETI method, see \cite{NME1620320604,MR2285892,MR2023942,MR2449391,MR2793906}. 
An alternative approach dubbed Multi-Trace formalism, leading to different solvers, 
was introduced a few years ago \cite{MR3273156,MR2856427,MR2927645,HJLP13_504,8601085,MR3069956,CH14_556,CHJ12_781}, 
providing other boundary integral formulations adapted to multi-domain geometrical 
configurations. Multi-trace boundary integral formulations seem well adapted to block 
preconditioners for domain decomposition but still very little is known about associated 
iterative global solvers. To our knowledge, the only contributions in this direction are
\cite{HJLP13_504,dolean:hal-00949024}.

A precise knowledge of the spectral structure of the equation under consideration
is most of the time a key ingredient for devising efficient domain decomposition
strategies. The main purpose of the present contribution is thus to provide results
concerning the spectrum of local multi-trace operators. In particular we describe the spectrum
through an explicit formula in the case where the propagation medium admits uniform
characteristic material parameters, which yields the location of
the accumulation points of the spectrum of local multi-trace operators in the general case
of piecewise varying effective wave numbers. In addition, the present contribution provides a proof
of the well posedness of the local multi-trace formulations in the case of non-trivial
relaxation parameters, which was mentioned in \cite[Rem.2]{HJLP13_504} as an open problem.

The present contribution is organized as follows. In the first section we describe the scattering 
problem and the geometrical configurations under consideration. In Section \ref{TraceSpaces}
and \ref{SummaryPotentialTheory} we introduce notations related to traces and integral 
operators, and recall well established results of classical potential theory. In Section 
\ref{DerivationLocalMTF} we recall the derivation of (relaxed) local multi-trace formulations 
as presented in \cite{HJLP13_504}, and prove uniqueness of the solutions to this formulation 
for general values of the relaxation parameter. Section \ref{Examples} presents detailed calculus
achieved in the case of two particular elementary geometrical configurations. 
Section \ref{SpectrumGalOp} is dedicated to the study of the spectrum of the local multi-trace
operator, and Section \ref{Numerics} will present numerical results confirming the theory.

\section{Setting of the problem}
In this section, we will mainly introduce notations, and write the problem 
under consideration, starting from a  precise description of the geometrical configurations 
we wish to examine. First of all, we consider a partition $\R^{d}:=\cup_{j=0}^{n}
\overline{\Omega}_{j}$ where the $\Omega_{j}$'s are Lipschitz domain, where $\R^{d} =\R,\R^{2}$ or $\R^{3}$.
We assume that each $\Omega_{j}$ is bounded except $\Omega_{0}$. Changing the numbering of sub-domains 
if necessary, we may assume without loss of generality that each $\Omega_{j}$ is connected.
We shall refer to the boundary  of each sub-domain by $\Gamma_{j}:=\partial\Omega_{j}$, and also set 
$\Gamma_{j,k}:=\Gamma_{j}\cap \Gamma_{k} = \partial\Omega_{j}\cap\partial\Omega_{k}$ 
to refer to interfaces. The union of all interfaces will be denoted 
$$
\Sigma := \mathop{\cup}_{j=0}^{n}\Gamma_{j} = \mathop{\cup}_{0\leq j < k\leq n}\Gamma_{j,k}\;\;.
$$
This set will be referred to as the skeleton of the partition. We make a further 
strong geometric hypothesis, assuming that the sub-domain partition under consideration 
does not involve any junction point, so that each $\Gamma_{j,k}$ is a closed Lipschitz manifold 
without boundary,
\begin{equation}\label{NoJunctionPoint}
\partial\Gamma_{j,k} = \emptyset\quad\forall j,k=0\dots n.
\end{equation}
\begin{pspicture}(0,5)(0,0)
\put(5,2.5){
\psellipse(2,0)(4.5,2.25)
\pscircle(0,0){1}
\pscircle(4,0){1}
\put(-0.25,-0.1){$\Omega_{1}$}
\put(1.75,-0.1){$\Omega_{2}$}
\put(3.75,-0.1){$\Omega_{3}$}
\put(7,-0.1){$\Omega_{0}$}

\put(5.25,-1.5){\textcolor{white}{\rule{20pt}{10pt}}}
\put(-0.25,0.8){\textcolor{white}{\rule{25pt}{10pt}}}
\put(3.6,0.8){  \textcolor{white}{\rule{22.5pt}{10pt}}}

\put(5.5,-1.5){$\Gamma_{0,2}$}
\put(-.15,0.9){$\Gamma_{1,2}$}
\put(3.65,0.9){ $\Gamma_{2,3}$}

}
\end{pspicture}
\paragraph{Sobolev spaces} We need to introduce a few usual notations related to standard 
Sobolev spaces. If $\omega\subset \R^{d}$ is any Lipschitz domain, we set 
$\mH^{1}(\omega):=\{v\in\mL^{2}(\omega),\,\nabla v\in \mL^{2}(\omega)\}$ equipped 
with the norm $\Vert v\Vert_{\mH^{1}(\omega)}^{2}:=\Vert v\Vert_{\mL^{2}(\omega)}^{2} + 
\Vert \nabla v\Vert_{\mL^{2}(\omega)}^{2}$, and $\mH^{1}(\Delta,\omega):=\{v\in
\mH^{1}(\omega),\;\Delta v\in \mL^{2}(\omega)\}$ equipped with 
$\Vert v\Vert_{\mH^{1}(\Delta,\omega)}^{2}:=\Vert v\Vert_{\mH^{1}(\omega)}^{2}+ 
\Vert \Delta v\Vert_{\mL^{2}(\omega)}^{2}$. In addition, if $\mH(\omega)$ refers to any 
of the above mentioned spaces, then  $\mH_{\loc}(\overline{\omega})$ will refer to the space 
of functions $v$ such that $\psi v\in \mH(\omega)$ for any $\psi\in \mathscr{C}^{\infty}_{\mrm{comp}}
(\R^{d}):=\{ \varphi\in \mathscr{C}^{\infty}(\R^{d}),\;\mrm{supp}(\varphi)\;\textrm{bounded}\}$.

For any Lipschitz open set $\omega\subset\R^{d}$, we shall refer to the space 
of Dirichlet traces $\mH^{1/2}(\partial\omega) := \{ v\vert_{\partial\omega}, \; 
v\in\mH^{1}(\omega)\}$ equipped with the norm $\Vert v\Vert_{\mH^{1/2}(\partial\omega)} := 
\mrm{min}\{ \Vert u\Vert_{\mH^{1}(\omega)},\; u\vert_{\partial\omega} = v\}$. The space of 
Neumann traces $\mH^{-1/2}(\partial\omega)$ will be defined as the dual to 
$\mH^{1/2}(\partial\omega)$ equipped with the corresponding canonical dual norm
$\Vert p\Vert_{\mH^{-1/2}(\partial\omega)} := \sup_{v\in\mH^{1/2}(\partial\omega)\setminus\{0\}}
\vert\langle p,v\rangle\vert/\Vert v\Vert_{\mH^{1/2}(\partial\omega)}$.

\paragraph{Transmission problem} We will consider a very standard wave scattering problem 
(so-called transmission problem), imposing Hemholtz equation in each sub-domain, as well as 
transmission conditions across interfaces:
find $u\in \mH^{1}_{\loc}(\R^{d})$ such that 
\begin{equation}\label{InitPb}
\left\{\begin{array}{l}
-\Delta u -\kappa_{j}^{2}u = 0\quad \textrm{in}\;\;\Omega_{j}\quad\forall j=0\dots n\\[5pt]
u - u_{\inc}\;\; \textrm{is outgoing in}\;\Omega_{0}\\\quad\\
u\vert_{\Gamma_{j}} - u\vert_{\Gamma_{k}} = 0\\[5pt]
\partial_{n_{j}}u\vert_{\Gamma_{j}} + \partial_{n_{k}}u\vert_{\Gamma_{k}} = 0
\quad\textrm{on}\;\;\Gamma_{j,k} = \Gamma_{j}\cap \Gamma_{k},\;\;\forall j,k=0\dots n\\[5pt]
\end{array}\right.
\end{equation}
In the equation above $u_{\inc}\in \mH^{1}_{\loc}(\R^{d})$ is a known source term 
of the  problem satisfying $-\Delta u_{\inc} -\kappa_{0}^{2} u_{\inc} = 0$ in $\R^{d}$.
In addition, we assume that $\kappa_{j}>0$ for all $j=0\dots n$.
Problem (\ref{InitPb}) is known to admit a unique solution, see \cite{VonPetersdorff} 
for example. The outgoing condition in (\ref{InitPb}) refers to Sommerfeld's radiation
condition, see \cite{MR841971}. A function $v\in\mH^{1}_{\loc}(\Delta,\Omega_{0})$ will be
said $\kappa$-outgoing radiating if $\lim_{r\to\infty}\int_{\partial\mrm{B}_{r}}\vert\partial_{r}v-
\imath\kappa v\vert^{2}d\sigma = 0$ where $\mrm{B}_{r}$ refers to the ball of radius $r$
and center $\mathbf{0}$, and $\partial_{r}$ is the partialderivative with respect to
the radial variable $r = \vert\bfx\vert$.

\paragraph{Trace operators} As this problem involves transmission conditions, and since we are interested 
in boundary integral formulations of it, we need to introduce suitable trace 
operators. According to \cite[Thm. 2.6.8 and Thm 2.7.7]{SauterSchwab}, every 
sub-domain $\Omega_{j}$ gives rise to continuous boundary trace operators 
$\gamma^{j}_{\dir}:\mH^{1}_{\loc}(\overline{\Omega}_{j})\to \mH^{1/2}(\partial\Omega_{j})$ and 
$\gamma^{j}_{\neu}:\mH^{1}_{\loc}(\Delta, \overline{\Omega}_{j})\to \mH^{-1/2}(\partial\Omega_{j})$
(so-called Dirichlet and Neumann traces) uniquely defined by
$$
\gamma_{\dir}^{j}(\varphi) := \varphi\vert_{\partial\Omega_{j}}\quad\textrm{and}\quad
\gamma_{\neu}^{j}(\varphi) := \bfn_{j}\cdot\nabla \varphi\vert_{\partial\Omega_{j}}
\quad\quad\forall\varphi\in \mathscr{C}^{\infty}(\overline{\Omega}_{j}).
$$
In the definition above $\bfn_{j}$ refers to the vector field normal to $\partial\Omega_{j}$ 
pointing toward the exterior of $\Omega_{j}$. We will also need a notation to refer to an
operator gathering both traces in a single array 
$$
\gamma^{j}(u) := (\gamma^{j}_{\dir}(v),
\gamma^{j}_{\neu}(v))
$$ 
We shall also refer to $\gamma^{j}_{\dir,c},\gamma^{j}_{\neu,c}$ defined in the same manner as 
$\gamma^{j}_{\dir},\gamma^{j}_{\neu}$ with traces taken from the exterior of $\Omega_{j}$. In addition, 
we set $\gamma^{j}_{c}(v) := (\gamma^{j}_{\dir,c}(v),\gamma^{j}_{\neu,c}(v))$. We will refer 
to mean values and jumps to these trace operators, setting 
$$
\{\gamma^{j}(u)\} := \frac{1}{2}\big(\,\gamma^{j}(u) + \gamma^{j}_{c}(u)\,\big)
\quad\textrm{and}\quad 
\lbr\gamma^{j}(u)\rbr := \gamma^{j}(u) - \gamma^{j}_{c}(u)
$$

\section{Trace spaces}\label{TraceSpaces}

We want to recast Problem \eqref{InitPb} into variational boundary integral equations 
set in trace spaces adapted to the present multi-sub-domain context. The most simple space 
we can introduce consists in the \emph{multi-trace space} \cite[Sect.~2.1]{MR3069956} i.e. 
the cartesian product of local traces:
\begin{equation}\label{mtsdef}
\begin{array}{l}  
  \mbH(\Sigma)\;:=\;\mbH(\Gamma_{0})
  \times\cdots \times\mbH(\Gamma_{n})\\[5pt]
  \textrm{where}\quad 
  \mbH(\Gamma_{j}) := 
  \mH^{+\frac{1}{2}}(\Gamma_{j})\times
  \mH^{-\frac{1}{2}}(\Gamma_{j})\;.
\end{array}
\end{equation}
 We endow each $\mbH(\Gamma_{j})$ with the norm $\left\|
(v,q)\right\|_{\mbH(\Gamma_{j})}:= (\Vert v\Vert_{\mH^{1/2}(\Gamma_{j})}^{2} + \Vert
q\Vert_{\mH^{-1/2}(\Gamma_{j})}^{2})^{1/2}$, and equip $\mbH(\Sigma)$ with
the norm naturally associated with the cartesian product
$$
\left\|\ctr{u}\right\|_{\mbH(\Sigma)} := \big( \;\Vert \ctr{u}_{0}
\Vert_{\mbH(\Gamma_{0})}^{2} +\cdots +
\Vert \ctr{u}_{n}\Vert_{\mbH(\Gamma_{n})}^{2} \;\big)^{\frac{1}{2}}
$$
for $\ctr{u} = (\ctr{u}_{0},\dots, \ctr{u}_{n})\in
\mbH(\Sigma)$\footnote[3]{Functions in Dirichlet trace spaces like
  $\mH^{1/2}(\partial\Omega_{j})$ will be denoted by $u,v,w$, whereas we
  use $p,q,r$ for Neumann traces. Small fraktur font symbols $\ctr u$, $\ctr v$,
  $\ctr w$ are reserved for Cauchy traces, with an integer subscript indicating
  restriction to a sub-domain boundary.}. In the sequel we shall repeatedly refer 
to the continuous operator $\gamma:\Pi_{j=0}^{n}\mH^{1}_{\loc}(\Delta,\overline{\Omega}_{j})$
$\to \mbH(\Sigma)$ defined by $\gamma(u):=(\gamma_{0}(u),\dots,\gamma_{n}(u))$, where 
$\Pi_{j=0}^{n}\mH^{1}_{\loc}(\Delta,\overline{\Omega}_{j})$ should be understood as the 
set of $u\in\mL^{2}_{\loc}(\R^{d})$ such that $u\vert_{\Omega_{j}}\in \mH^{1}_{\loc}(\Delta,
\overline{\Omega}_{j})$ for all $j$. We also need a bilinear duality pairing
for $\mbH(\Gamma_{j})$ and $\mbH(\Sigma)$;   writing $\langle\,,\,\rangle_{\Gamma_{j}}$ 
for the duality pairing between $\mH^{1/2}(\Gamma_{j})$ and $\mH^{-1/2}(\Gamma_{j})$, 
we opt for the skew-symmetric bilinear form
\begin{equation}  \label{BilForm1}
\begin{array}{l}
\dsp{   \llbr\ctr{u},\ctr{v}\rrbr := 
  \sum_{j=0}^{n}\;\lbr\ctr u_{j},\ctr v_{j}\rbr_{\Gamma_{j}} }
  \quad\quad\textrm{where}\\[5pt]
  \dsp{ \Big\lbr
    \left(\begin{array}{c}
        u_{j}\\p_{j}
      \end{array}\right),
    \left(\begin{array}{c}
        v_{j}\\q_{j}
      \end{array}\right)\Big\rbr_{\Gamma_{j}}\;:=\;
  \langle u_{j},q_{j}\rangle_{\Gamma_{j}}-\langle v_{j},p_{j}\rangle_{\Gamma_{j}}\;. }\\
\end{array}
\end{equation}
Next, as in \cite[Sect.~2.2]{MR3069956}, \cite[Sect.~3.1]{CHJ12_781}, we introduce
the so-called single-trace space that consists in collections of traces complying 
with transmission conditions. This space can be defined by 
\begin{equation}\label{DefSingleTraceSpace}
\mbX(\Sigma) := \mrm{clos}_{\,\mbH(\Sigma)} \{\;\gamma(u)=(\gamma^{j}(u))_{j=0}^{n}\;
\vert\; u\in \mathscr{C}^{\infty}(\R^{d})\;\}
\end{equation}
where the symbol $\mrm{clos}_{\,\mbH(\Sigma)}$ refers to the closure with respect to
the norm on $\mbH(\Sigma)$.
By construction, this is a closed subspace of $\mbH(\Sigma)$. 
Note also that a function $u\in \mH^{1}_{\loc}(\Delta,\overline{\Omega}_{0})\times
\dots\times\mH^{1}_{\loc}(\Delta,\overline{\Omega}_{n})$ satisfies the transmission conditions
of (\ref{InitPb}), if and only if $\gamma(u) = (\gamma^{j}(u))_{j=0}^{n}\in\mbX(\Sigma)$. In
particular, if $u\in \mH^{1}_{\loc}(\Delta,\R^{d})$ then $\gamma(u) =(\gamma^{j}(u))_{j=0}^{n}\in\mbX(\Gamma)$. 
In the sequel, we will use this space to enforce transmission conditions. The single-trace 
space admits a simple weak characterization, see  \cite[Prop.2.1]{MR3069956}.

\begin{lemma}\label{CaracSingleTrace}\quad\\
For any $\ctr{u}\in\mbH(\Sigma)$ we have $\ctr{u}\in \mbX(\Sigma)\iff \llbr \ctr{u},\ctr{v}\rrbr = 0\;\forall \ctr{v}\in \mbX(\Sigma)$.
\end{lemma}

\section{Summary of potential theory}\label{SummaryPotentialTheory}

In this paragraph, we shall remind the reader of well established results
concerning the integral representation of solutions to homogeneous Helmholtz 
equation in Lipshitz domains. A detailed proof of the statements contained in the 
present paragraph can be found for example in \cite[Chap.3]{SauterSchwab}.\\

\noindent
Let $\kernel_{\kappa}(\bfx)$ refer to the outgoing Green's kernel associated to the Helmholtz 
operator $-\Delta-\kappa^{2}$. For example $\kernel_{\kappa}(\bfx) = \exp(i\kappa\vert\bfx\vert)
/(4\pi\vert\bfx\vert)$ in $\R^{3}$. For each sub-domain $\Omega_{j}$, any $(v,q)\in 
\mbH(\Gamma_{j})$ and any $\bfx\in\R^{d}\setminus\Gamma_{j}$, define the potential 
operator
\begin{equation}\label{PotOp}
\Green_{\kappa}^{j}(v,q)(\bfx) := \int_{\Gamma_{j}}q(\bfy)\,\kernel_{\kappa}(\bfx-\bfy) + 
v(\bfy)\,\bfn_{j}(\bfy)\cdot(\nabla\kernel_{\kappa})(\bfx-\bfy)d\sigma(\bfy)\;.
\end{equation}
The operator $\Green_{\kappa}^{j}$ maps continuously 
$\mbH(\Gamma_{j})$ into $\mH^{1}(\Delta,\overline{\Omega}_{j})\times 
\mH^{1}(\Delta,\R^{d}\setminus\Omega_{j})$, see \cite[Thm 3.1.16]{SauterSchwab}. In particular 
$\Green_{\kappa}^{j}$ can be applied to a pair of traces of the form $\ctr{v} = \gamma^{j}(v)$. 
This potential operator can be used to write a representation formula for solution to 
homogeneous Helmholtz equations, see \cite[Thm 3.1.6]{SauterSchwab}.

\begin{proposition}\label{ReprThm}\quad\\
Let $u\in \mH^{1}_{\loc}(\overline{\Omega}_{j})$ satisfy $\Delta u+\kappa_{j}^{2}u = 0$ 
in $\Omega_{j}$. Assume in addition that $u$ is $\kappa_{0}$-outgoing if $j=0$. We have the 
representation formula 
\begin{equation}\label{ReprFormula}
\Green_{\kappa_{j}}^{j}(\gamma^{j}(u))(\bfx)\;=\;
\left\{\begin{array}{ll}
u(\bfx) & \textrm{for}\;\;\bfx\in \Omega_{j}\\[5pt]
0       & \textrm{for}\;\;\bfx\in \R^{d}\setminus\overline{\Omega}_{j}\;.
\end{array}\right.
\end{equation}
Similarly, if $v\in \mH^{1}_{\loc}(\R^{d}\setminus\Omega_{j})$ satisfies $\Delta v+\kappa_{j}^{2}v= 0$ 
in $\R^{d}\setminus \overline{\Omega}_{j}$, and is $\kappa_{j}$-outgoing radiating in the case $j\neq 0$, 
then we have $\Green_{\kappa_{j}}^{j}(\gamma_{c}^{j}(v))(\bfx) = -v(\bfx)$ for $\bfx\in \R^{d}\setminus
\overline{\Omega}_{j}$, and $\Green_{\kappa_{j}}^{j}(\gamma_{c}^{j}(v))(\bfx) = 0$ for $\bfx\in \Omega_{j}$.
\end{proposition}

\noindent 
The potential operator (\ref{PotOp}) also satisfies a remarkable identity, known as jump
formula, describing the behavior of $\Green_{\kappa_{j}}^{j}(\ctr{v})(\bfx)$ as $\bfx$ crosses 
$\Gamma_{j}$, namely
\begin{equation}\label{JumpFormulas}
\lbr \gamma^{j}\rbr\cdot \Green_{\kappa_{j}}^{j}(\ctr{v}) = \ctr{v}\quad \quad\forall \ctr{v}\in \mbH(\Gamma_{j})
\end{equation}
which also writes  $\lbr \gamma^{j}\rbr\cdot \Green_{\kappa_{j}}^{j} = \mrm{Id}$.
We refer the reader to \cite[Thm.3.3.1]{SauterSchwab} (the jump formulas are more commonly 
written in the form of four equations in the literature). Proposition \ref{ReprThm} shows that, 
if $u$ is solution to a homogeneous Helmholtz equation in $\Omega_{j}$ (and is outgoing if $j=0$) 
then $\gamma^{j}\cdot \Green_{\kappa_{j}}^{j}(\gamma^{j}(u)) = \gamma^{j}(u)$. This actually turns 
out to be a characterization of traces of solutions to homogeneous Helmholtz equation. 

\begin{proposition}\label{ProjCald}\quad\\
Define $\mathcal{C}_{\kappa}^{\mrint}(\Omega_{j}):=
\{\gamma^{j}(u)\;\vert\; u\in\mH^{1}_{\loc}(\overline{\Omega}_{j})\;,\;\Delta u +\kappa^{2}u=0\;\textrm{in}\;\Omega_{j},\;
\textrm{$u$ outgoing if $j=0$}\;\}$. Then $\gamma^{j}\cdot\Green_{\kappa}^{j}:\mbH(\Gamma_{j})\to 
\mbH(\Gamma_{j})$ is a continuous projector, called Calder\'on projector interior to $\Omega_{j}$, 
whose range coincides with $\mathcal{C}_{\kappa}^{\mrint}(\Omega_{j})$ i.e. for any $\ctr{v}\in\mbH(\Gamma_{j})$
$$
\gamma^{j}\cdot\Green_{\kappa}^{j}(\ctr{v}) = \ctr{v}\quad\iff
\quad \ctr{v} \in \mathcal{C}_{\kappa}^{\mrint}(\Omega_{j})\;.
$$
Similarly, defining $\mathcal{C}_{\kappa}^{\mrout}(\Omega_{j}):=
\{\gamma^{j}_{c}(u)\;\vert\; u\in\mH^{1}_{\loc}(\overline{\Omega}_{j}),\;\Delta u +\kappa^{2}u=0\;\textrm{in}\;
\R^{d}\setminus\Omega_{j},\;\textrm{$u$ outgoing if } j \neq 0\;\}$, 
we have $\gamma^{j}\cdot\Green_{\kappa}^{j}(\ctr{v}) = 0$ if and only if
$\ctr{v} \in \mathcal{C}_{\kappa}^{\mrout}(\Omega_{j})$.
\end{proposition}
\noindent 
For a detailed proof, see  \cite[Prop.3.6.2]{SauterSchwab}.
We shall repeatedly use this characterization as a convenient way to express wave equations 
in the sub-domains $\Omega_{j}$. Here is another characterization of 
the space of Cauchy data which was established in \cite[Lemma 6.2]{MR3069956}.

\begin{lemma}\label{CaracCauchy}\quad\\
Consider any $j=0,\dots n$, and any $\kappa\in\mathbb{C}\setminus\{0\}$ such that 
$\Re e\{\kappa\}\geq 0$, $\Im m\{\kappa\}\geq 0$. Then for any $\ctr{u}\in\mbH(\Gamma_{j})$ 
we have: $\ctr{u}\in \mathcal{C}_{\kappa}^{\mrint}(\Omega_{j})\iff \lbr\ctr{u},\ctr{v}\rbr_{\Gamma_{j}} = 0\;
\forall \ctr{v}\in \mathcal{C}_{\kappa}^{\mrint}(\Omega_{j})$. Similarly we have
$\ctr{u}\in \mathcal{C}_{\kappa}^{\mrout}(\Omega_{j})\iff \lbr\ctr{u},\ctr{v}\rbr_{\Gamma_{j}} = 0\;
\forall \ctr{v}\in \mathcal{C}_{\kappa}^{\mrout}(\Omega_{j})$. 
\end{lemma}
\quad\\[-10pt]
The results that we have stated above hold for any $j=0\dots n$. We also set 
$\mathcal{C}^{\alpha}(\Sigma):=\mathcal{C}^{\alpha}_{\kappa_{0}}(\Omega_{0})\times\cdots
\times \mathcal{C}^{\alpha}_{\kappa_{n}}(\Omega_{n})$ for $\alpha = \mrint,\mrout$.
The notations just introduced allow a condensed reformulation of the 
well-posedness of (\ref{InitPb}), see  \cite[Prop.6.1]{MR3069956} for
a detailed discussion and proof.
\begin{lemma}\label{DirectSum}\quad\\
$\mbX(\Sigma)\oplus \mathcal{C}^{\mrint}(\Sigma) = \mbH(\Sigma)$.
\end{lemma}

\noindent
To handle Calder\'on projectors in a multi-sub-domain context, it is more comfortable 
to introduce global operators, so as to reduce notations. First we introduce the continuous 
operator $\mrA:\mbH(\Sigma)\to \mbH(\Sigma)$ defined by
\begin{equation}\label{DefGlobOp}
\begin{array}{l}
\llbr \mrA(\ctr{u}),\ctr{v}\rrbr := \sum_{j=0}^{n} \left\lbr \mrA^{j}(\ctr{u}_{j}), 
\ctr{v}_{j})\right\rbr_{\Gamma_{j}}\\[10pt]
\textrm{with}\quad\mrA^{j} := 2\,\{\gamma^{j}\}\cdot\Green^{j}_{\kappa_{j}}
\end{array}
\end{equation}
For all  $\ctr{u} = (\ctr{u}_{j})_{j=0}^{n}, \ctr{v} = (\ctr{v}_{j})_{j=0}^{n}\in\mbH(\Sigma)$. 
Observe that $(\mrm{Id}\pm\mrA)/2$ are projectors, according to Proposition \ref{ProjCald}, with 
$\mrm{ker}(\mrA-\mrm{Id}) = \mrm{range}(\mrA+\mrm{Id}) = \mathcal{C}^{\mrint}(\Sigma)$ and 
$\mrm{ker}(\mrA+\mrm{Id}) = \mrm{range}(\mrA-\mrm{Id}) = \mathcal{C}^{\mrout}(\Sigma)$. 
With this notation, a direct consequence of Proposition \ref{ProjCald} is that 
$(\mrA^{j})^{2} = (2\{\gamma^{j}\}\cdot\Green_{\kappa}^{j} )^{2} = \mrm{Id}$. This is summarized 
in the next lemma.

\begin{lemma}\label{CalderonIdentity}\quad\\
$(\mrA)^{2} = \mrm{Id}$.
\end{lemma}

\begin{remark}\label{RemarkTwoDomains}
In the case of two sub-domains $\R^{d} = \overline{\Omega}_{0}\cup \overline{\Omega}_{1}$ 
and one interface $\Sigma = \Gamma_{0} = \Gamma_{1}$, with $\kappa = \kappa_{0} = \kappa_{1}$,
there is a remarkable identity relating $\mrA^{0}$ to $\mrA^{1}$.  
Indeed, in this situation, the only difference in the definition of 
$\mrm{G}^{0}_{\kappa}$ and $\mrm{G}^{1}_{\kappa}$ comes from $\bfn_{0} = -\bfn_{1}$.
In particular, denoting 
$$
Q := 
\left\lbr\begin{array}{cc}
1 & 0 \\ 0 & -1
\end{array}\right\rbr
$$
we have $\mrm{G}^{0}_{\kappa}(\ctru) = -\mrm{G}^{1}_{\kappa}(Q\cdot\ctru)\;
\forall \ctru\in \mH^{1/2}(\Gamma_{0})\times \mH^{-1/2}(\Gamma_{0})$. Note that 
$\{\gamma^{1}\} = Q\cdot\{\gamma^{0}\}$, so multiplying the previous equality by 
$\{\gamma^{1}\}$ yields $Q\cdot\mrA^{0} = -\mrA^{1}\cdot Q$.\hfill $\Box$

\end{remark}

\noindent 
Lemma \ref{CalderonIdentity} above shows directly that $\mrA$ is invertible. It also satisfies a generalized 
Garding inequality. The next result is proved for example in \cite[\S 4]{VonPetersdorff}.

\begin{proposition}\label{GardingA}\quad\\
Define  $\Theta:\mbH(\Sigma)\to \mbH(\Sigma)$ by 
$\Theta(\ctr{u}) := (Q(\ctr{u}_{0}), Q(\ctr{u}_{1}),\dots, Q(\ctr{u}_{n}))$ 
for all $\ctr{u} = (\ctr{u}_{0},\ctr{u}_{1},\dots, \ctr{u}_{n})$, 
where $Q(u,p) := (u,-p)$ for any $(u,p)\in \mH^{1/2}(\Gamma_{j})\times 
\mH^{-1/2}(\Gamma_{j})$ and any $j$. There exists a compact operator 
$\mrm{K}:\mbH(\Sigma)\to \mbH(\Sigma)$ and a constant $C>0$ such that
$$
\Re e\{-\llb \mrA\ctru,\Theta(\overline{\ctru})\rrb + 
\llb \mrm{K}\ctru,\overline{\ctru}\rrb\,\}\geq C\,
\Vert \ctru\Vert_{\mbH(\Sigma)}^{2}\quad \quad\forall 
\ctru\in\mbH(\Sigma).
$$
\end{proposition}

\section{Local multi-trace formulation}\label{DerivationLocalMTF}
In this section we would like to recall the derivation of the local multi-trace formulation
introduced in \cite{HJLP13_504}, and provide detailed analysis for it. The formulation considered 
in  \cite{MR2927645} is a particular case of the formulations introduced in \cite{HJLP13_504}. 
But only an analysis of the formulation of \cite{MR2927645} has been provided so far. In addition,
we should underline our hypothesis (\ref{NoJunctionPoint}) that discards any junction point in the 
sub-domain partition, while \cite{MR2927645,HJLP13_504} also considered geometrical configurations 
with junction points.

\quad\\
A key ingredient of local multi-trace theory is an operator yielding a characterization of transmission conditions 
of (\ref{InitPb}). Considering $\ctru = (u_{k},p_{k})_{k=0}^{n}$ and $\ctrv = (v_{j},q_{j})_{j=0}^{n}$,
we define the transmission operator $\Pi:\mbH(\Sigma)\to\mbH(\Sigma)$ by
\begin{equation}\label{DefTransOp}
\ctr{v} = \Pi(\ctr{u})\;\;\iff\;\;
\left\{\begin{array}{l}
v_{j} = +\,u_{k}\\[5pt]
q_{j} = -\,p_{k}
\end{array}\right.
\quad \textrm{on}\;\;\Gamma_{j,k}\quad\forall j,k = 0\dots n, 
\;\;j\neq k.
\end{equation}
\quad\\
Clearly, for any function $u\in \mH^{1}_{\loc}(\Delta,\R^{d})$ we have $\gamma(u) := 
(\gamma^{j}(u))_{j=0}^{n}\in \mbH(\Sigma)$ with $\gamma(u) = \Pi(\gamma(u))$. 
Conversely, considering any function $u\in \mL^{2}_{\loc}(\R^{d})$ such that 
$u\vert_{\Omega_{j}}\in \mH^{1}_{\loc}(\Delta,\overline{\Omega}_{j})$ for all $j=0\dots n$, 
then $\gamma(u)\in \mbH(\Sigma)$ is well defined, and if $\gamma(u) = \Pi(\gamma(u))$ 
then $u\in \mH^{1}_{\loc}(\Delta,\R^{d})$. Routine calculus shows the following remarkable identities
\begin{equation}\label{RemId}
\Pi^{2} = \mrm{Id}\;,\quad \overline{\Pi(\ctr{v})} =  \Pi(\overline{\ctr{v}})\quad\textrm{and}\quad \llbr \Pi(\ctr{u}),\ctr{v}\rrbr = 
\llbr \Pi(\ctr{v}),\ctr{u}\rrbr\quad \forall \ctr{u},\ctr{v}\in\mbH(\Sigma). 
\end{equation}
As is readily checked, the operator $\Pi$ maps continuously $\mbH(\Sigma)$ onto $\mbH(\Sigma)$
\textit{under Assumption (\ref{NoJunctionPoint}) that each $\Gamma_{j,k}$ is a Lipschitz manifold 
without boundary}. Elementary arguments on trace spaces show that,
for any $\ctr{u}\in \mbH(\Sigma)$, we have $\ctr{u}\in\mbX(\Sigma)\iff 
\ctr{u}  =\Pi(\ctr{u})$. Since $\Pi^{2} = \mrm{Id}$, this can be simply rewritten 
in the following manner.

\begin{lemma}\label{AltCarac}\quad\\
$\mrm{range}(\Pi+\mrm{Id}) = \mbX(\Sigma)$.
\end{lemma}

\quad\\[-10pt]
Now consider $\ctr{u} = (\gamma^{j}(u))_{j=0}^{n}$ the traces of the unique solution $u$ 
to Problem (\ref{InitPb}). The homogeneous wave equation satisfied by $u$ in each sub-domain
can be reformulated by means of Calder\'on projectors $(\mrm{A}-\mrm{Id})(\ctr{u}-\ctr{u}_{\inc}) = 0$. 
Choose any relaxation parameter $\alpha\in \C\setminus\{0\}$ and add $\alpha\,(\mrm{Id} - \Pi)\ctru = 0$ 
to this equation, which is consistent since $\ctru = \gamma(u)$ must satisfy the transmission 
conditions of (\ref{InitPb}). Denoting $\ctr{f}:= (\mrm{A}-\mrm{Id})\ctr{u}_{\inc}$, what precedes implies 
\begin{equation}\label{LocMltTr}
\begin{array}{l}
\left\{\begin{array}{l}
\ctr{u}\in \mbH(\Sigma)\quad\textrm{and}\\[5pt]
\llbr\;(\mrA-\Pi_{\alpha}) \ctr{u},\ctr{v}\;\rrbr = \llbr \ctr{f},\ctr{v}\rrbr\quad\forall \ctr{v}\in \mbH(\Sigma).
\end{array}\right.\\[20pt]
\quad \textrm{where}\quad \Pi_{\alpha}:= (1-\alpha)\mrm{Id} + \alpha \Pi
\end{array}
\end{equation}
Observe that the operator of the formulation above can also be rewritten in the form of a 
convex combination $\mrA-\Pi_{\alpha} =(1-\alpha)(\mrA-\mrm{Id}) + \alpha(\mrA-\Pi)$. Existence and uniqueness 
of the solution to this formulation has already been established 
only in the case $\alpha = 1$, see  \cite[Thm.9 and Thm.11]{MR2927645}. For all other 
values of $\alpha$, well-posedness of this formulation was an open problem so far, as mentioned 
in \cite[Rem.2]{HJLP13_504}. Below we prove uniqueness of the solution to (\ref{LocMltTr}) for \textit{any} 
$\alpha\in \C\setminus\{0\}$.

\begin{proposition}\label{Uniqueness}\quad\\
$\mrm{Ker}(\mrA-\Pi_{\alpha}) = \{0\}
\quad\forall \alpha\in\C\setminus\{0\}$.
\end{proposition}

\noindent \textbf{Proof:}

Take $\ctr{u} = (\ctr{u}_{j})_{j=0}^{n}$ satisfying $(\mrA-\mrm{Id})\ctr{u} + \alpha(\mrm{Id}-\Pi)\ctr{u} = 0$. Thus
we have  $\ctr{w} := \alpha (\Pi-\mrm{Id})\ctr{u} = (\mrA-\mrm{Id})\ctr{u}\in 
\mrm{range}(\mrA-\mrm{Id})\cap\mrm{range}(\Pi-\mrm{Id})$. Denote $\ctr{w}_{j}$ the component of $\ctr{w}$ associated 
to $\partial\Omega_{j}$ so that $\ctr{w} = (\ctr{w}_{0},\cdots,\ctr{w}_{n})$, and set 
$\psi_{j}(\bfx) := \mrm{G}^{j}_{\kappa_{j}}(\ctr{w}_{j})(\bfx)$. We have $(\mrA+\mrm{Id})\ctr{w} = (\mrA^{2} - \mrm{Id})\ctr{u} = 0$, 
which can be rewritten $\gamma^{j}\cdot \mrm{G}^{j}_{\kappa_{j}}(\ctr{w}_{j}) =\gamma^{j}(\psi_{j}) = 0$ for all 
$j=0\dots n$. Since we also have, by construction, $-\Delta\psi_{j}-\kappa_{j}^{2}\psi_{j} = 0$ in $\Omega_{j}$, 
we conclude that $\psi_{j} = 0$ in $\Omega_{j}$, and $\ctr{w}_{j} = \lbr\gamma^{j}\rbr \cdot \mrm{G}^{j}_{\kappa_{j}}(\ctr{w}_{j}) 
= \lbr\gamma^{j}(\psi_{j})\rbr = -\gamma^{j}_{c}(\psi_{j})$. If we can prove that $\psi_{j} = 0$ in $\R^{d}\setminus
\overline{\Omega}_{j}$ for each $j$, this will show that $\ctr{w} = 0$. 

\quad\\
Now observe that, since $\ctr{w}\in \mrm{range}(\Pi-\mrm{Id})$ we have $\Pi(\ctr{w}) + \ctr{w} = 0$ i.e 
$\Pi(\ctr{w}) = - \ctr{w}$. As a consequence the functions $\psi_{j}$ satisfy an homogeneous problem 
with "anti-transmission conditions"
\begin{equation}\label{DualBVP}
\left\{\begin{array}{l}
-\Delta \psi_{j} -\kappa_{j}^{2}\psi_{j} = 0\quad\textrm{in}\;\;\R^{d}\setminus\overline{\Omega}_{j}
\quad \forall j=0\dots n\;,\\[5pt]
\psi_{j}\;\textrm{is outgoing (with respect to $\kappa_{j}$) for $j\neq 0$}\\[5pt]
\gamma^{j}_{\dir,c}(\psi_{j}) + \gamma^{k}_{\dir,c}(\psi_{k}) = 0\quad\textrm{and}\\[5pt]
\gamma^{j}_{\neu,c}(\psi_{j}) - \gamma^{k}_{\neu,c}(\psi_{k}) = 0\quad 
\textrm{on}\;\;\Gamma_{j}\cap\Gamma_{k}\quad \forall j,k.
\end{array}\right.
\end{equation}
\quad\\
Since $\Pi^{2} = \mrm{Id}$ and $\Pi(\ctr{w}) = -\ctr{w}$, we have 
$2\ctr{w}  = \ctr{w}-\Pi(\ctr{w})$, and $\overline{\ctr{w}}+\Pi(\overline{\ctr{w}}) = 0$.
From this and (\ref{RemId}), we obtain $2 \llbr \ctr{w},\overline{\ctr{w}}\rrbr = \llbr \ctr{w}-\Pi(\ctr{w}),
\overline{\ctr{w}}\rrbr = -\llbr \overline{\ctr{w}}+\Pi(\overline{\ctr{w}}), \ctr{w}\rrbr = 0$. This can 
be rewritten 
\begin{equation}\label{SumEqualZero}
\sum_{j=0}^{n}\Im m\Big\{\,\int_{\Gamma_{j}}\gamma_{\dir,c}^{j}(\psi_{j})
\gamma_{\neu,c}^{j}(\overline{\psi}_{j})d\sigma\,\Big\}=0.
\end{equation}
Take $r>0$ sufficiently large to garantee that $\R^{d}\setminus\Omega_{0}\subset \mrB_{r}$ where 
$\mrB_{r}\subset \R^{d}$ refers to the ball centered at $0$ with radius $r$. Since $-\Delta\psi_{j} - 
\kappa_{j}^{2}\psi_{j} = 0$ in $\R^{d}\setminus\overline{\Omega}_{j}$, applying Green's formula 
in each $\mrB_{r}\setminus\overline{\Omega}_{j}$ yields 
$$
\begin{array}{rl}
   \dsp{  \int_{\partial\mrB_{r}}\psi_{j}\partial_{r}\overline{\psi}_{j} d\sigma }
&  \dsp{= \int_{\mrB_{r}\setminus\Omega_{j}}\vert\nabla \psi_{j}\vert^{2} - \kappa_{j}^{2}\vert\psi_{j}\vert^{2}d\bfx
 +  \int_{\partial\Omega_{j}}\gamma^{j}_{\dir,c}(\psi_{j})\gamma^{j}_{\neu,c}(\overline{\psi}_{j}) d\sigma }\\[10pt]

 0 
&  \dsp{= \int_{\mrB_{r}\setminus\Omega_{0}}\vert\nabla \psi_{0}\vert^{2} - \kappa_{0}^{2}\vert\psi_{0}\vert^{2}d\bfx
 +  \int_{\partial\Omega_{0}}\gamma^{0}_{\dir,c}(\psi_{0})\gamma^{0}_{\neu,c}(\overline{\psi}_{0}) d\sigma }

\end{array}
$$
In these equations, "$\partial_{r}$" refers to the radial derivative. Take the imaginary part 
of the identities above, and sum over $j=0\dots n$, taking account of (\ref{SumEqualZero}). This 
leads to 
$$
\sum_{j=0}^{n}\Im m\Big\{\int_{\partial\mrB_{r}}\psi_{j}\partial_{r}\overline{\psi}_{j} d\sigma\Big\} =
\Im m\Big\{ \sum_{j=0}^{n}\int_{\partial\Omega_{j}}\gamma^{j}_{\dir,c}(\psi_{j})
\gamma^{j}_{\neu,c}(\overline{\psi}_{j}) d\sigma\Big\} = \frac{1}{2\imath}\llb \ctr{w},\overline{\ctr{w}} \rrb = 0
$$
where $\imath$ refers to the imaginary unit. By construction, the functions $\psi_{j}$ are outgoing radiating, 
so that $0 = \lim_{r\to \infty}\int_{\partial\mrB_{r}}\vert\partial_{r}\psi_{j} - \imath\kappa_{j}\psi_{j}\vert^{2}d\sigma = 0 $. 
As a consequence, we finally obtain

$$
\begin{array}{l}
\dsp{ \sum_{j=0}^{n}\frac{1}{\kappa_{j}}\int_{\partial\mrB_{r}}
\vert\partial_{r}\psi_{j}\vert^{2} + \kappa_{j}^{2}\vert\psi_{j}\vert^{2}d\sigma}\\[10pt]

\hspace{1.5cm} \dsp{ = \sum_{j=1}^{n}\frac{1}{\kappa_{j}}\int_{\partial\mrB_{r}}
\vert\partial_{r}\psi_{j}-i\kappa_{j}\psi_{j}\vert^{2}d\sigma + 2\sum_{j=1}^{n}
\Im m\Big\{\int_{\partial\mrB_{r}} \psi_{j}\partial_{r}\overline{\psi}_{j}d\sigma\Big\} }\\[10pt]

\hspace{1.5cm} \dsp{ = \sum_{j=1}^{n}\frac{1}{\kappa_{j}}\int_{\partial\mrB_{r}}
\vert\partial_{r}\psi_{j}-i\kappa_{j}\psi_{j}\vert^{2}d\sigma\;
\mathop{\longrightarrow}_{r\to \infty}\; 0 }

\end{array}
$$
This shows in particular that $\lim_{r\to \infty}\int_{\partial\mrB_{r}}\vert\psi_{j}\vert^{2}d\sigma = 0$ 
for all $j=1\dots n$. As a consequence we can apply Rellich's lemma, see \cite[Lemma 3.11]{MR700400}, which 
implies that $\psi_{j} = 0$ in the unbounded connected component of each $\R^{d}\setminus\overline{\Omega}_{j}$.

Let us show that $\psi_{j}$ also vanishes in bounded connected components of $\R^{d}\setminus\overline{\Omega}_{j}$.
Take an arbitrary $j$, and let $\mathcal{O}$ be a bounded connected component of $\R^{d}\setminus\overline{\Omega}_{j}$. 
We have $\partial\mathcal{O} = \partial\Omega_{j}\cap \partial\Omega_{k} = \Gamma_{j,k}$ for some $k=0\dots n,k\neq j$.
Let $\mathcal{O}'$ be the unbounded connected component of $\R^{d}\setminus\overline{\Omega}_{k}$. Then we have 
$\Omega_{j}\subset \mathcal{O}'$, and $\partial\mathcal{O}' = \partial\mathcal{O} = \partial\Omega_{j}\cap 
\partial\Omega_{k}$. Since $\psi_{k} = 0$ in $\mathcal{O}'$, we have $\gamma^{j}_{\dir,c}(\psi_{j})\vert_{\Gamma_{j,k}} = 
- \gamma^{k}_{\dir,c}(\psi_{k}) = 0$ and $\gamma^{j}_{\neu,c}(\psi_{j})\vert_{\Gamma_{j,k}} = \gamma^{k}_{\neu,c}(\psi_{k}) = 0$,
according to the transmission conditions of (\ref{DualBVP}). Finally we have $-\Delta\psi_{j}-\kappa_{j}^{2}\psi_{j} = 0$
in $\mathcal{O}$ with $\gamma^{j}_{c}(\psi_{j}) = 0$ on $\partial\mathcal{O}$. We conclude by unique continuation 
principle (see \cite[\S.4.3]{MR841971}) that $\psi_{j} = 0$ in $\mathcal{O}$. We have just proved that 
$$
\psi_{j} = 0\quad\textrm{in}\;\;\R^{d}\setminus\overline{\Omega}_{j}\quad\forall j=0\dots n.
$$
Finally we have, on the one hand, $\ctr{w} = (\mrA - \mrm{Id})\ctr{u} = 0$, 
so $\ctr{u} = (\mrA + \mrm{Id})\ctr{u}/2\in \mrm{range}(\mrA + \mrm{Id}) = \mathcal{C}^{\mrint}(\Sigma)$, 
and on the other hand $\ctr{w} = \alpha(\Pi - \mrm{Id})\ctr{u} = 0$, so $\ctr{u} = (\Pi + \mrm{Id})\ctr{u}/2
\in \mrm{range}(\Pi + \mrm{Id}) =\mbX(\Sigma)$. So we conclude that $\ctr{u}\in \mathcal{C}^{\mrint}(\Sigma)
\cap \mbX(\Sigma) = \{0\}$ according to Lemma \ref{DirectSum}. \hfill $\Box$

\section{Examples}\label{Examples}

Before going further into the analysis of the local multi-trace formulation (\ref{LocMltTr}), 
we dedicate this section to deriving and studying it in ultra simplified situations 
where all calculations can be conducted quasi-explicitly. Here we will systematically
consider the case where all wave numbers are equal
\begin{equation}\label{AllWNEqual}
\kappa_{0} = \kappa_{1} = \dots  = \kappa_{n}.
\end{equation}
This assumption will allow substantial simplifications. Another purpose of
the present section is to determine the spectrum of the multi-trace operator 
in these simplified situations.
% Preliminary calculations, as well as numerical 
%results, have already been presented in \cite{???} for the two cases considered below. 
%Here we adopt a more systematic presentation though. 

\subsection{Two domain configuration}
We start by considering the case where the space is partitioned in two domains only. 
This simple case was already considered in \cite[\S 3.1]{MR2927645}, but here we are going to 
formulate additional comments. In this case $\Sigma = \Gamma_{0} = \Gamma_{1} = \Gamma_{0,1}$. 
We want to represent the operator $(1-\alpha)( \mrA-\mrm{Id}) + \alpha(\mrA-\Pi):\mbH(\Gamma_{0,1})\times 
\mbH(\Gamma_{0,1}) \to \mbH(\Gamma_{0,1})\times \mbH(\Gamma_{0,1})$ in a matrix form. First of 
all, note that the operator $\Pi$ admits the following expression,
\begin{equation}\label{ExplicitPi}
\Pi\,(\left\lbr\begin{array}{c}
\ctr{v}_{0}\\\ctr{v}_{1}
\end{array}\right\rbr) = 
\left\lbr\begin{array}{cc}
0 & \mrQ \\ \mrQ & 0
\end{array}\right\rbr \cdot
\left\lbr\begin{array}{c}
\ctr{v}_{0}\\\ctr{v}_{1}
\end{array}\right\rbr\quad\quad 
\textrm{with}\quad\quad Q:=
\left\lbr\begin{array}{cc}
1 & 0 \\ 0 & -1
\end{array}\right\rbr.
\end{equation}
Hence, denoting
$\ctr{u} =(\ctr{u}_{0},\ctr{u}_{1})$, we have 
$$
(\mrA-\alpha\Pi)\ctr{u}\;=
\left\lbr\begin{array}{cc}
\mrA^{0} & -\alpha\mrQ \\[5pt]
-\alpha\mrQ & \mrA^{1}
\end{array}\right\rbr \cdot
\left\lbr\begin{array}{c}
\ctr{u}_{0}\\[5pt] 
\ctr{u}_{1}
\end{array}\right\rbr.
$$
To determine the spectrum of the (relaxed) multi-trace operator $(1-\alpha)( \mrA-\mrm{Id}) + \alpha(\mrA-\Pi)$, 
it suffices to determine the spectrum of  $\mrA - \alpha\Pi$. If we compute the square of this operator, taking 
account of (\ref{AllWNEqual}), we obtain $(\mrA-\alpha\Pi)^{2} = \mrA^{2} +\alpha^{2}\Pi^{2} -\alpha 
(\Pi\mrA + \mrA\Pi) = (1+\alpha^{2})\mrm{Id} - \alpha(\Pi\mrA + \mrA\Pi)$. In this case, a 
direct calculus shows that $\mrQ\mrA^{0} = -\mrA^{1}\mrQ$. As a consequence, an explicit 
calculus yields
$$
\begin{array}{ll}
\left\lbr\begin{array}{cc}
\mrA^{0} & 0 \\ 0 & \mrA^{1}
\end{array}\right\rbr \cdot
\left\lbr\begin{array}{cc}
0 & \mrQ \\ \mrQ & 0 
\end{array}\right\rbr 

& = \left\lbr\begin{array}{cc}
0 & \mrA^{0}\mrQ \\ \mrA^{1}\mrQ & 0
\end{array}\right\rbr \\[15pt]

& = \left\lbr\begin{array}{cc}
0 & -\mrQ\mrA^{1} \\ -\mrQ\mrA^{0} & 0
\end{array}\right\rbr = 

- \left\lbr\begin{array}{cc}
0 & \mrQ \\ \mrQ & 0 
\end{array}\right\rbr \cdot
\left\lbr\begin{array}{cc}
\mrA^{0} & 0 \\ 0 & \mrA^{1}
\end{array}\right\rbr 

\end{array}
$$
From this we conclude that $\Pi\mrA + \mrA\Pi = 0$, which finally yields 
$(\mrA-\alpha\Pi)^{2} = (1+\alpha^{2})\,\mrm{Id}$. This expression, together with
the spectral mapping theorem ( see \cite[Thm.10.28]{MR1157815} for example),
provides an explicit characterization of the spectrum in the case where all wave-numbers 
are equal\footnote[4]{In the remaining of this article, the square root 
of complex numbers shall be defined by $\sqrt{\rho\exp(i\theta)} = \sqrt{\rho} \exp(i\theta/2)$
for $\theta\in (-\pi,\pi\rbr$.  }
$$
\mathfrak{S}(\mrA-\Pi_{\alpha}) \subset
\{-1+\alpha+\sqrt{1+\alpha^{2}}, -1+\alpha-\sqrt{1+\alpha^{2}}\}.
$$

\subsection{Three domain configuration}

Now we consider a partition in three domains $\R^{d} = \overline{\Omega}_{0}\cup 
\overline{\Omega}_{1} \cup \overline{\Omega}_{2}$. This situation is pictured below. 
\quad\\
\begin{pspicture}(0,5.5)(0,0)
\put(7,3){

\pspolygon(-2,-2)(2,-2)(2,2)(-2,2)
\pscircle(0,0){1}

\put(-0.1,-0.1){$\Omega_{1}$}
\put(1,-1.25){$\Omega_{0}$}
\put(3,-0.1){$\Omega_{2}$}

\put(-2.2,0){\textcolor{white}{\rule{10pt}{17.5pt}}}
\put(-2.2,0.2){$\Gamma_{0,2}$}
\put(0.1,0.5){\textcolor{white}{\rule{17.5pt}{15pt}}}
\put(0.2,0.9){$\Gamma_{0,1}$}

}
\end{pspicture}
\noindent
As the definition of the transmission operator is given by a formula 
written on each interface, let us decompose traces on each sub-domain according 
to  interfaces. Considering the decomposition $\Gamma_{0} = \Gamma_{0,1}\cup 
\Gamma_{0,2}$, any trace $\ctr{v}\in \mbH(\Gamma_{0})$ induces 
an element $\mrm{R}^{0}_{1}(\ctr{v})\in \mbH(\Gamma_{0,1})$ defined by 
$\mrm{R}^{0}_{1}(\ctr{v}) = \ctr{v}\vert_{\Gamma_{0,1}}$. We may define 
$\mrm{R}^{0}_{2}(\ctr{v})\in \mbH(\Gamma_{0,2})$ similarly. This establishes a natural isomorphism 
$(\mrm{R}^{0}_{1},\mrm{R}^{0}_{2}):\mbH(\Gamma_{0})\to\mbH(\Gamma_{0,1})\times \mbH(\Gamma_{0,2})$. 
The adjoint of those maps are extension operators i.e. $(\mrm{R}^{0}_{j})^{*}(\ctr{v}) = \ctr{v}
\cdot 1_{\Gamma_{0,j}}\in \mbH(\Gamma_{0})$ for any $\ctr{v}\in \mbH(\Gamma_{0,j})$. With these maps, the operator $\mrA^{0}$ 
induces a $2\times 2$ matrix denoted $\lbr \mrA^{0}\rbr$ with integral operator entries
$$
\lbr\mrA^{0}\rbr = 
\left\lbr\begin{array}{cc}
\mrA^{0}_{1,1} & \mrA^{0}_{1,2}\\[10pt]
\mrA^{0}_{2,1} & \mrA^{0}_{2,2}
\end{array}\right\rbr\quad\quad 
\textrm{with} \quad \mrA^{0}_{j,k}:=
\mrm{R}^{0}_{j}\cdot\mrA^{0}\cdot(\mrm{R}^{0}_{k})^{*}
$$
Plugging this decomposition into the definition of the operator $\mrm{A}$ yields  
a $4\times 4$ matrix of integral operators acting on tuples $(\ctr{u}_{0,1},\ctr{u}_{0,2},
\ctr{u}_{1},\ctr{u}_{2})\in \mbH(\Gamma_{0,1})\times\mbH(\Gamma_{0,2})\times \mbH(\Gamma_{1})
\times \mbH(\Gamma_{2})$, and given by the following formula  
\begin{equation}\label{ExprExplOp}
(\mrA-\alpha\Pi)\ctr{u} =
\left\lbr\begin{array}{cccc}
\mrA^{0}_{1,1} & \mrA^{0}_{1,2} & -\alpha\mrQ & 0 \\[5pt]
\mrA^{0}_{2,1} & \mrA^{0}_{2,2} & 0 & -\alpha\mrQ \\[5pt]
-\alpha\mrQ & 0 & \mrA^{1} & 0 \\[5pt]
0 & -\alpha\mrQ & 0 & \mrA^{2} \\
\end{array}\right\rbr\cdot
\left\lbr\begin{array}{c}
\ctr{u}_{0,1} \\[5pt] \ctr{u}_{0,2} \\[5pt] \ctr{u}_{1} \\[5pt] \ctr{u}_{2} 
\end{array}\right\rbr
\end{equation}
As in the previous paragraph, let us compute the spectrum of the local multi-trace operator
$(1-\alpha)(\mrA-\mrm{Id}) + \alpha(\mrA-\Pi)$. Here again, 
it suffices to determine the spectrum of $(\mrA-\alpha\Pi)$. 
Once again we have $(\mrm{A}-\alpha\Pi)^{2} = \mrm{A}^{2} + \alpha^{2}\Pi^{2} - 
\alpha(\Pi\mrA+\mrA\Pi) = (1+\alpha^{2})\,\mrm{Id} - \alpha(\Pi\mrA+\mrA\Pi)$. Besides, 
taking account of (\ref{AllWNEqual}), a direct and thorough calculus shows that 
$\mrQ\cdot\mrA^{0}_{j,j}\cdot\mrQ = -\mrA^{j}$. So if we compute explicitly the expression 
of $\Pi\mrA+\mrA\Pi$ taking account of this identity, we obtain
$$
\Pi\mrA+\mrA\Pi = 
\left\lbr\begin{array}{cccc}
0 & 0 & 0 & \mrA^{0}_{1,2}\mrQ \\[5pt]
0 & 0 & \mrA^{0}_{2,1}\mrQ & 0 \\[5pt]
0 & \mrQ\mrA^{0}_{1,2} & 0 & 0 \\[5pt]
\mrQ\mrA^{0}_{2,1} & 0 & 0 & 0 \\
\end{array}\right\rbr
$$
This time we have $\Pi\mrA+\mrA\Pi\neq 0$. Let us compute the square 
of this operator. Since $\mrQ^{2} = \mrm{Id}$, we obtain
$$
(\Pi\mrA+\mrA\Pi)^{2} = 
\left\lbr\begin{array}{cccc}
\mrA^{0}_{1,2}\mrA^{0}_{2,1} & 0 & 0 & 0 \\[5pt]
0 & \mrA^{0}_{2,1}\mrA^{0}_{1,2} & 0 & 0 \\[5pt]
0 & 0 & \mrQ\mrA^{0}_{1,2}\mrA^{0}_{2,1}\mrQ & 0 \\[5pt]
0 & 0 & 0 & \mrQ\mrA^{0}_{2,1}\mrA^{0}_{1,2}\mrQ \\
\end{array}\right\rbr
$$
Now let us have a closer look at each of the operators 
$\mrA^{0}_{1,2}\mrA^{0}_{2,1}$ and $\mrA^{0}_{2,1}\mrA^{0}_{1,2}$.
Observe that $-\mrA^{0}_{j,k} = 2\,\mrQ\cdot\gamma^{j}\cdot\Green^{k}\cdot\mrQ$ for $j\neq k$.
According to the second equality in (\ref{ReprFormula}), 
we have $\gamma^{1}\cdot\Green^{2}\cdot\gamma^{2}\cdot\Green^{1} = 0$.
As a consequence, we obtain 
$$
\begin{array}{ll}
\mrA^{0}_{1,2}\mrA^{0}_{2,1} 
& = 4\,\big(\mrQ\cdot\gamma^{1}\cdot\Green^{2}\cdot\mrQ\big)\cdot 
    \big(\mrQ\cdot\gamma^{2}\cdot\Green^{1}\cdot\mrQ\big)\\[5pt]
& = 4\,\mrQ\cdot(\gamma^{1}\cdot\Green^{2}\cdot\gamma^{2}\cdot\Green^{1})\cdot\mrQ = 0
\end{array}
$$
We show in a similar manner that $\mrA^{0}_{2,1}\mrA^{0}_{1,2} $. To conclude we have 
$(\Pi\mrA + \mrA\Pi)^{2} = 0$. Such a nilpotent operator has a non-empty spectrum
(see \cite[Thm.10.13]{MR1157815}) that is reduced to $\{0\}$ according to
the spectral mapping theorem \cite[Thm.10.28]{MR1157815}, which implies of course
that $\mathfrak{S}(\Pi\mrA + \mrA\Pi) = \{0\}$. Finally we obtain 
the following spectrum, like in the previous paragraph, 
$$
\mathfrak{S}(\mrA-\Pi_{\alpha}) \subset
\big\{\;-1+\alpha+\sqrt{1+\alpha^{2}}, -1+\alpha-\sqrt{1+\alpha^{2}}\;\big\}.
$$

\section{Spectrum of the operator in a general configuration}\label{SpectrumGalOp}

For both examples of the previous two paragraphs, the spectrum of the local 
multi-trace operator only consisted in the two eigenvalues 
$-1+\alpha\pm\sqrt{1+\alpha^{2}}$ in the case where all wave 
numbers equal. Besides, during the calculations above, the geometry
of the interfaces never came into play. In the present section we will
show that these are actually general results that hold for any
number of sub-domain arranged arbitrarily, provided the geometry
does not involve any junction point.

\quad\\
To investigate this question in the general case, we need to introduce 
further notations. Recall that each boundary can be decomposed in the 
following manner $\Gamma_{j} = \cup_{k\neq j}\Gamma_{j,k}$. We will decompose 
traces accordingly. For a given pair $j,k$ with $j\neq k$ we define 
$$
\begin{array}{l}
\mrR^{j}_{k}:\mbH(\Gamma_{j})\to \mbH(\Gamma_{j,k})\\[5pt]
\mrR^{j}_{k}(\ctr{v}) := \ctr{v}\vert_{\Gamma_{j,k}}
\end{array}
$$
To reformulate the above definition fully explicitly, for $\ctr{v} = (v,q)\in \mbH(\Gamma_{j,k}) 
= \mH^{1/2}(\Gamma_{j,k})\times \mH^{-1/2}(\Gamma_{j,k})$ and $\ctr{u} = (u,p)\in \mbH(\Gamma_{j})
= \mH^{1/2}(\Gamma_{j})\times \mH^{-1/2}(\Gamma_{j})$, we have $\ctr{v} = \mrR^{j}_{k}(\ctr{u})$ if 
and only if $u\vert_{\Gamma_{j,k}} = v$ and $p\vert_{\Gamma_{j,k}} = q$. 
The adjoint of these restriction operators are given by
$$
\begin{array}{l}
(\mrR^{j}_{k})^{*}:\mbH(\Gamma_{j,k})\to \mbH(\Gamma_{j})\\[5pt]
(\mrR^{j}_{k})^{*}(\ctr{v}) := \ctr{v}\cdot\mrm{1}_{\Gamma_{j,k}}
\end{array}
$$
Decomposing traces with the embedding/restriction operators that we have 
just defined, each $\mrA^{j}:\mbH(\Gamma_{j})\to \mbH(\Gamma_{j})$ 
induces a matrix of integral operators denoted $\lbr \mrA^{j}\rbr$ with 
maximal size $(n-1)\times(n-1)$ given by 
\begin{equation}\label{DecompOp}
\lbr \mrA^{j}\rbr := 
\left\lbr\begin{array}{ccc}
\mrA^{j}_{0,0} & \cdots & \mrA^{j}_{0,n}\\
\vdots & & \vdots\\
\mrA^{j}_{n,0} & \cdots & \mrA^{j}_{n,n}
\end{array}\right\rbr\quad\quad
\mrA^{j}_{k,m}:= \mrR^{j}_{k}\cdot\mrA^{j}\cdot(\mrR^{j}_{m})^{*}
\end{equation}
In the notation above, it should be understood that the rows and the columns 
associated to indices $k$ such that $\Gamma_{j,k} = \partial\Omega_{j}\cap 
\partial\Omega_{k} = \emptyset$ must be omitted. The row/column associated to 
$k=j$ is to be omitted as well. For each sub-domain $\Omega_{j}$, we will need to 
consider the set of indices
$$
\mathscr{I}_{j}:=\{\;k\;\vert\; k\neq j, \;\partial\Omega_{j}\cap \partial\Omega_{k}
\neq \emptyset\}.
$$
Hence the matrix in (\ref{DecompOp}) is square with $\mrm{card}(\mathscr{I}_{j})$ rows.
If we plug the definition of the restriction/embedding operators $\mrR^{j}_{k}$ 
into the definition (\ref{DefGlobOp}) of the operator $\mrA^{j}$ and the potential 
operator (\ref{PotOp}), we obtain an explicit formula for $\mrA^{j}_{k,m}$ in the 
case where $k\neq m$ namely,
\begin{equation}\label{ExplicitFormulaExtraDiag}
\textcolor{white}{+}\begin{array}{l}
\dsp{ -\Big\llb\; \mrA^{j}_{k,m}
\Big(\begin{array}{c}
u \\ p \end{array}\Big),
\Big(\begin{array}{c}
v \\ q \end{array}\Big)
\;\Big\rrb_{\Gamma_{j,k}}\;=}  \\[10pt]  

\textcolor{white}{+}\dsp{ \int_{\Gamma_{j,k}} \int_{\Gamma_{j,m}}\Green_{\kappa}(\bfx-\bfy)\big\lbr\; 
(\bfn(\bfx)\times\nabla u(\bfx))\cdot(\bfn(\bfy)\times\nabla v(\bfy)) }\\

\hspace{4cm}\dsp{ -\kappa^{2}\bfn(\bfx)\cdot\bfn(\bfy)u(\bfx)v(\bfy) 
-p(\bfx)q(\bfy)\;\big\rbr\;d\sigma(\bfy )d\sigma(\bfx)}\\[10pt]

+ \dsp{ \int_{\Gamma_{j,k}} \int_{\Gamma_{j,m}}(\nabla\Green_{\kappa})(\bfx-\bfy)\cdot\big\lbr\;
q(\bfy)u(\bfx)\bfn(\bfx) - p(\bfx)v(\bfy)\bfn(\bfy) \;\big\rbr\;d\sigma(\bfy )d\sigma(\bfx)   }
\end{array}
\end{equation}
for any $(u,p)\in\mbH(\Gamma_{j,m})$ and $(v,q)\in \mbH(\Gamma_{j,k})$. Note that the expression 
above does not hold for $k=m$. Expression (\ref{ExplicitFormulaExtraDiag}) clearly shows that 
the operators $\mrA^{j}_{k,m}$ are compact for $k\neq m$ since it only involves smooth kernels. 
The following lemma yields several remarkable identities 
satisfied by the elements of the decomposition (\ref{DecompOp}).
\begin{lemma}\label{ElementaryId}\quad\\
For any $j=0\dots n$, and any $k,l,m\in\mathscr{I}_{j}$ with $k\neq m$ and $k\neq l$  we have 
\begin{itemize}
\item[i)]   $(\mrA^{j}_{k,k})^{2} = \mrm{Id}$, 
\item[ii)]  $\mrQ\cdot\mrA^{j}_{k,k} = -\mrA^{k}_{j,j}\cdot\mrQ$ if $\kappa_{j} = \kappa_{k}$
\item[iii)] $\mrA^{j}_{l,k}\cdot\mrA^{j}_{k,m} = 0$,

\end{itemize}
\end{lemma}

\noindent \textbf{Proof:}

Pick an arbitrary $j = 0\dots n$ that will be fixed until the end of the proof.
Take a arbitrary $k,m\in \mathscr{I}_{j}$ with $k\neq m$. Let $\mathcal{O}\subset\R^{d}$
be the maximal open set satisfying $\partial\mathcal{O} = \Gamma_{j,m}$ and $\Omega_{j}\cap
\mathcal{O} = \emptyset$. Take an arbitrary $\ctr{v}\in
\mbH(\Gamma_{j,m})$, and denote $\tilde{\ctr{v}} :=(\mrR^{j}_{m})^{*}\ctr{v}$. Observe 
that $\Green^{j}_{\kappa}(\tilde{\ctr{v}})\in \mH^{2}_{\loc}(\R^{d}\setminus\mathcal{O})$. Since 
$\Gamma_{j,k}\cap\Gamma_{j,m} = \emptyset$, this implies in particular that
$\Green^{j}_{\kappa}(\tilde{\ctr{v}})$ does not admit any jump across $\Gamma_{j,k}$. 
As a consequence, we have 
$$
\mrA^{j}_{k,m}(\ctr{v}) = 2\,\mrR^{j}_{k}\cdot \gamma^{j}\cdot\Green^{j}(\tilde{\ctr{v}}) = 
2\,\mrR^{j}_{k}\cdot \gamma^{j}_{c}\cdot\Green^{j}(\tilde{\ctr{v}})
$$
Now observe that, if $\ctr{w}\in \mathcal{C}^{\mrout}_{\kappa}(\Omega_{j})$, then 
$(\mrR^{j}_{k})^{*}\mrR^{j}_{k}(\ctr{w})\in \mathcal{C}^{\mrout}_{\kappa}(\Omega_{j})$ for any $k \in\mathscr{I}_{j}$. 
Taking $\ctr{w} = \gamma^{j}_{c}\cdot\Green^{j}(\tilde{\ctr{v}})$, we see that 
$\tilde{\ctr{w}} =(\mrR^{j}_{k})^{*}\mrA^{j}_{k,m}(\ctr{v})\in \mathcal{C}^{\mrout}_{\kappa}(\Omega_{j})$. 
According to the integral representation Theorem \ref{ReprThm}, this implies 
$\Green^{j}_{\kappa}(\tilde{\ctr{w}})(\bfx) = 0$ for $\bfx\in\Omega_{j}$. In particular we have 
$$
2\,\mrR^{j}_{l}\gamma^{j}\Green^{j}_{\kappa}(\tilde{\ctr{w}}) = \mrA^{j}_{l,k}\cdot\mrA^{j}_{k,m}(\ctr{v}) = 0.
$$
This establishes \textit{iii)}. Now we know that $(\mrA^{j})^{2} = \mrm{Id}$. This implies 
that, for any $k \in\mathscr{I}_{j}$, we have $\sum_{m\in\mathscr{I}_{j}} \mrA^{j}_{k,m}\cdot \mrA^{j}_{m,k}= \mrm{Id}$.
But according to \textit{iii)}, all the terms of this sum vanish, except for $k=m$ which establishes \textit{i)}.

\quad\\
To prove \textit{ii)}, observe that $\mrQ\cdot\mrR^{j}_{k}\cdot\{\gamma^{j}\} = \mrR^{k}_{j}\cdot\{\gamma^{k}\}$. 
On the other hand, the explicit expression of potential operators given by (\ref{PotOp}) shows that 
$\Green^{j}_{\kappa}\cdot(\mrR^{j}_{k})^{*} = -\Green^{k}_{\kappa}\cdot(\mrR^{k}_{j})^{*}\cdot \mrQ$ in the case where 
$\kappa_{j} = \kappa_{k}$. Combining these two identities we obtain $\mrQ\cdot\mrA^{j}_{k,k} = 
2\,\mrQ\cdot\mrR^{j}_{k}\cdot\{\gamma^{j}\}\cdot\Green^{j}_{\kappa}\cdot(\mrR^{j}_{k})^{*} =
-2\,\mrR^{k}_{j}\cdot\{\gamma^{k}\} \cdot \Green^{k}_{\kappa}\cdot(\mrR^{k}_{j})^{*}\cdot \mrQ = -\mrA^{k}_{j,j}\cdot \mrQ$.

\hfill $\Box$

\quad\\
Next we need to introduce an operator involving only the diagonal blocks of the matrix 
representing $\mrA$ in the decomposition (\ref{DecompOp}), without any term coupling 
different interfaces. Define $\mrm{D}:\mbH(\Sigma)\to\mbH(\Sigma)$  by the 
explicit formula 
\begin{equation}\label{DefOpD}
\begin{array}{l}
\dsp{ \llb \mrm{D}(\ctr{u}),\ctr{v}\rrb := \sum_{j=0}^{n}\sum_{k\in \mathscr{I}_{j}}
\;\llb\; \mrA^{j}_{k,k}(\ctr{u}_{j,k}), \ctr{v}_{j,k}\;\rrb_{\Gamma_{j,k}}\quad\quad 
\forall \ctr{u},\ctr{v}\in \mbH(\Sigma) }\\[20pt]

\textrm{where}\quad \ctr{u}_{j,k} := \mrR^{j}_{k}(\ctr{u}_{j}),\quad \ctr{v}_{j,k} := \mrR^{j}_{k}(\ctr{v}_{j})
\end{array}
\end{equation}
\begin{lemma}\label{IdentityDiagOp}\quad\\
We have $(\mrm{D})^{2} = \mrm{Id}$. Moreover, if all wave numbers equal i.e. 
$\kappa_{0} = \kappa_{1} = \dots =\kappa_{n}$, we have  $\Pi\cdot\mrm{D} + 
\mrm{D}\cdot\Pi = 0$.
\end{lemma}
\noindent\textbf{Proof:}

According to Lemma \ref{ElementaryId}, we already know that $(\mrA^{j}_{k,k})^{2} = \mrm{Id}$.  
Pick an arbitrary pair of traces $\ctr{u},\ctr{v}\in\mbH(\Sigma)$, and set $\ctr{u}_{j,k}:=\mrR^{j}_{k}(\ctr{u}_{j})$ and 
$\ctr{v}_{j,k}:=\mrR^{j}_{k}(\ctr{v}_{j})$ for all $j$, and all $k\in\mathscr{I}_{j}$.  Applying Formula (\ref{DefOpD}) 
twice yields 
$$
\begin{array}{ll}
\llb \mrm{D}^{2}(\ctr{u}),\ctr{v}\rrb 
& \dsp{ := \sum_{j=0}^{n}\sum_{k\in \mathscr{I}_{j}}\llb (\mrA^{j}_{k,k})^{2}\ctr{u}_{j,k}, \ctr{v}_{j,k}\rrb_{\Gamma_{j,k}}  }\\
& \dsp{ = \sum_{j=0}^{n}\sum_{k\in \mathscr{I}_{j}}\llb \ctr{u}_{j,k}, \ctr{v}_{j,k}\rrb_{\Gamma_{j,k}} = \sum_{j=0}^{n} 
\llb \ctr{u}_{j}, \ctr{v}_{j}\rrb_{\Gamma_{j}} = \llb \ctr{u}, \ctr{v}\rrb. }
\end{array}
$$
which establishes that $\mrm{D}^{2} = \mrm{Id}$. To establish the second statement, 
let us first point out that the definition (\ref{DefTransOp}) of the operator $\Pi$ 
can be rewritten $\llb \Pi(\ctr{u}),\ctr{v} \rrb = \sum_{j=0}^{n}\sum_{k\in\mathscr{I}_{j}}
\llb \mrQ(\ctr{u}_{k,j}), \ctr{v}_{j,k}\rrb_{\Gamma_{j,k}}$. Combining this with the 
definition of the operator $\mrm{D}$, and using Property \textit{iii)} of Lemma 
\ref{ElementaryId}, we obtain
$$
\begin{array}{ll}
\llb \Pi\cdot\mrm{D}(\ctr{u}),\ctr{v}\rrb 
& \dsp{ := \sum_{j=0}^{n}\sum_{k\in \mathscr{I}_{j}}\llb \mrQ\,\mrA^{k}_{j,j}\ctr{u}_{k,j}, \ctr{v}_{j,k}\rrb_{\Gamma_{j,k}}  }\\
& \dsp{ := \sum_{j=0}^{n}\sum_{k\in \mathscr{I}_{j}}-\llb \mrA^{j}_{k,k}\,\mrQ(\ctr{u}_{k,j}), \ctr{v}_{j,k}\rrb_{\Gamma_{j,k}} 
  =-\llb \mrm{D}\cdot\Pi(\ctr{u}), \ctr{v}\rrb 
  }\\
\end{array}
$$
\hfill $\Box$

\begin{lemma}\quad\\
The operator $\mrm{T}:=\mrA - \mrm{D}$ is compact and satisfies $\mrm{T}^{2} = 0$.
\end{lemma}
\noindent\textbf{Proof:}

The operator $\mrT$ only involves terms $\mrA^{j}_{k,m}$
with $k\neq m$. Since $\Gamma_{j,k}\cap \Gamma_{j,m} = \emptyset$, these operators defined 
by (\ref{ExplicitFormulaExtraDiag}) only involve smooth kernels. So each $\mrA^{j}_{k,m}, k\neq m$ 
is compact, and $\mrT$ is compact itself. Let us compute explicitly the expression of 
$\mrm{T}^{2}$. Pick arbitrary $\ctr{u},\ctr{v}\in\mbH(\Sigma)$, and set 
$\ctr{u}_{j,k}:=\mrR^{j}_{k}(\ctru_{j})$ and $\ctr{v}_{j,k}:=\mrR^{j}_{k}(\ctrv_{j})$ 
for all $j$, and all $k\in\mathscr{I}_{j}$. We have
$$
\begin{array}{ll}
\llb \mrm{T}(\ctr{u}),\ctr{v}\rrb &    \dsp{ := \sum_{j=0}^{n}\;\sum_{\substack{k,m\in \mathscr{I}_{j}\\ k\neq m  }  }
\llb \mrA^{j}_{k,m}(\ctr{u}_{j,m}), \ctr{v}_{j,k}\rrb_{\Gamma_{j,k}}   }\\[10pt]

\llb \mrm{T}^{2}(\ctr{u}),\ctr{v}\rrb & \dsp{ := \sum_{j=0}^{n}\;
\sum_{\substack{k,l,m\in \mathscr{I}_{j}\\ k\neq l, l\neq m   }  }
\llb \mrA^{j}_{k,l}\mrA^{j}_{l,m}(\ctr{u}_{j,m}), \ctr{v}_{j,k}\rrb_{\Gamma_{j,k}}}

\end{array}
$$
To conclude it remains to apply property \textit{ii)} of Lemma \ref{ElementaryId}. Since, 
for each term $\mrA^{j}_{k,l}\mrA^{j}_{l,m}$ we have $k\neq l$ and $m\neq l$, the whole 
sum vanishes, and we have $\llb \mrm{T}^{2}(\ctr{u}),\ctr{v}\rrb = 0$. \hfill $\Box$

\newpage \noindent 
Note that it is a direct consequence of the above lemma and Proposition \ref{GardingA} that 
$\mrm{D}$ also satisfies a Garding inequality. In addition, since $\Pi^{2} = \mrm{Id}$, a recurrence 
argument combined with the previous lemma readily leads to the following corollary. Note that, 
remarkably, this result holds \textit{without} any assumption on the wave numbers.

\begin{corollary}\label{Nilpotent1}\quad\\
We have $(\Pi\mrm{T} + \mrm{T}\Pi)^{k} = (\Pi\mrm{T})^{k} + (\mrm{T}\Pi)^{k},\quad\forall k\geq 0$.
\end{corollary}

\noindent 
Now let us formulate a few elementary and useful remarks concerning the geometrical arrangement 
of the interfaces. Let $\Upsilon = \{0,1,\dots,n\}$, and say that two indices
$j,k$ are adjacent if $\partial\Omega_{j}\cap \partial\Omega_{k}\neq \emptyset$.
This adjacency relation endow $\Upsilon$ with a graph structure. We order the elements of $\Upsilon$, 
writing $j\prec k$ if $\Omega_{j}$ is included in a bounded connected component of $\R^{d}\setminus
\overline{\Omega}_{k}$. This induces an tree struture on $\Upsilon$. In particular $\Upsilon$ is 
a tree and does not admit any (simple cycle), see for example \cite[Chap.16]{MR1855295}. 
And it does not admit chain with a length larger than $n$. The picture below provides 
an example of such a tree structure. 

\begin{pspicture}(0,7)(0,0)
\put(2,3.75){
\psellipse(2,0)(4,2.25)
\psellipse(0.5,0)(2,1)
\pscircle(4,0){1}
\pscircle(1,0){0.75}
\put(-0.5,-0.1){$\Omega_{1}$}
\put(1.75,-1.5){$\Omega_{2}$}
\put(3.75,-0.1){$\Omega_{3}$}
\put(0.8,-0.1){$\Omega_{4}$}
\put(-1.5,-2){$\Omega_{0}$} }

\put(9,4){

\psline(2,1.5)(2,0)
\psline(2,0)(1,-1)
\psline(1,-1)(1,-2.5)
\psline(2,0)(3,-1)

\pscircle[fillstyle=solid,fillcolor=white](2,1.5){0.35}
\pscircle[fillstyle=solid,fillcolor=white](2,0){0.35}
\pscircle[fillstyle=solid,fillcolor=white](1,-1){0.35}
\pscircle[fillstyle=solid,fillcolor=white](3,-1){0.35}
\pscircle[fillstyle=solid,fillcolor=white](1,-2.5){0.35}

\put(1.9,1.35){0}
\put(1.9,-.15){2}
\put(.9,-1.15){1}
\put(2.9,-1.15){3}
\put(.9,-2.65){4}
}

\end{pspicture}

\begin{proposition}\label{Nilpotent2}\quad\\
In the case where all wave numbers equal $\kappa_{0} = \kappa_{1} = \dots = \kappa_{n}$,
we have $(\Pi\mrm{T})^{n} = (\mrm{T}\Pi)^{n} = 0$ where $n$ is the number of sub-domains.
\end{proposition}

\noindent \textbf{Proof:}

First of all note that $(\mrm{T}\Pi)^{n} = \Pi (\Pi\mrm{T})^{n}\Pi$, so we only need to prove 
the result for $(\Pi\mrm{T})^{n}$. We start by simply writing down the explicit expression 
of the operator $(\Pi\mrm{T})^{n}$. For any  $\ctr{u},\ctr{v}\in\mbH(\Sigma)$, setting 
$\ctr{u}_{j,k}:=\mrR^{j}_{k}(\ctru_{j})$ and $\ctr{v}_{j,k}:=\mrR^{j}_{k}(\ctrv_{j})$, we have 
$$
\begin{array}{ll}
\llb \Pi\mrm{T}(\ctr{u}),\ctr{v}\rrb &    \dsp{ := \sum_{i_{0}=0}^{n}\;
\sum_{i_{1}\in \mathscr{I}_{i_{0}} }\sum_{\substack{  i_{2}\in \mathscr{I}_{i_{1}}\\ i_{2}\neq i_{0}  } }
\llb \mrQ\mrA^{i_{1}}_{i_{0},i_{2}}(\ctr{u}_{i_{1},i_{2}}), \ctr{v}_{i_{0},i_{1}}\rrb_{\Gamma_{i_{0},i_{1}}}   }\\[10pt]

\llb (\Pi\mrm{T})^{2}(\ctr{u}),\ctr{v}\rrb &    \dsp{ := \sum_{i_{0}=0}^{n}\;
\sum_{i_{1}\in \mathscr{I}_{i_{0}} }\sum_{\substack{  i_{2}\in \mathscr{I}_{i_{1}}\\ i_{2}\neq i_{0}  } }
\sum_{\substack{  i_{3}\in \mathscr{I}_{i_{2}}\\ i_{3}\neq i_{1}  } }
\llb \mrQ\mrA^{i_{1}}_{i_{0},i_{2}} \mrQ\mrA^{i_{2}}_{i_{1},i_{3}}(\ctr{u}_{i_{2},i_{3}}), \ctr{v}_{i_{0},i_{1}}\rrb_{\Gamma_{i_{0},i_{1}}}   }

\end{array}
$$
In the expression above we have $i_{2}\neq i_{0}$ and $i_{3}\neq i_{1}$ due to the very definition
of the operator $\mrm{T}$. Applying recursively the formulas derived above for $\Pi\mrm{T}$ 
finally leads to the following explicit expression for $(\Pi\mrm{T})^{n}$,

\noindent 
\begin{equation}\label{ExpressionCompliquee}
\begin{array}{ll}
\llb (\Pi\mrm{T})^{n}(\ctr{u}),\ctr{v}\rrb &    \dsp{ := \sum_{i_{0}=0}^{n}\;
\sum_{i_{1}\in \mathscr{I}_{i_{0}} }\sum_{\substack{  i_{2}\in \mathscr{I}_{i_{1}}\\ i_{2}\neq i_{0}  } }
\sum_{\substack{  i_{3}\in \mathscr{I}_{i_{2}}  \\ i_{3}\neq i_{1}        } }   
\sum_{\substack{  i_{4}\in \mathscr{I}_{i_{3}}  \\ i_{3}\neq i_{2}    } } \cdots\cdots 
\sum_{\substack{  i_{n+1}\in \mathscr{I}_{i_{n}}\\ i_{n+1}\neq i_{n-1}  } }    }\\[30pt]

& \hspace{2.5cm} \dsp{\llb \mrQ\mrA^{i_{1}}_{i_{0},i_{2}} \mrQ\mrA^{i_{2}}_{i_{1},i_{3}} \cdots
\mrQ\mrA^{i_{n}}_{i_{n-1},i_{n+1}}
(\ctr{u}_{i_{n},i_{n+1}}), \ctr{v}_{i_{0},i_{1}}\rrb_{\Gamma_{i_{0},i_{1}}}   }

\end{array}
\end{equation}
The sum in the expression above is taken over all the sequence of indices $i_{0},i_{1},
\dots,i_{n+1}$ satisfying the constraints $i_{k}\in \mathscr{I}_{i_{k-1}}$ (which implies 
in particular that $i_{k}\neq i_{k-1}$) and $i_{k+1}\neq i_{k-1}$ for all $k=1\dots n$. 
Each $i_{k}$ is the index of the sub-domain $\Omega_{i_{k}}$, and $\Omega_{i_{k+1}}$ is adjacent to 
$\Omega_{i_{k}}$ since $i_{k+1}\in\mathscr{I}_{i_{k}}$, hence those sequences
$i_{0},i_{1},\dots,i_{n+1}$ are actually chains of length exactly $n+1$ of the tree 
$\Upsilon$. But since $\Upsilon$ only admits $n$ elements and is a tree, it does not contain 
any such chain. This implies that the summation in (\ref{ExpressionCompliquee}) is taken over an 
empty set. Hence $\llb (\Pi\mrm{T})^{n}(\ctr{u}),\ctr{v}\rrb = 0$, and since 
$\ctr{u},\ctr{v}$ were chosen arbitrarily, this finally implies $(\Pi\mrm{T})^{n} = 0$.
\hfill $\Box$

\quad\\
The spectrum of the local multi-trace operator can now easily be deduced
from what precedes, in the case where all wave numbers equal. Note that the next result 
states equality and not just inclusion. 

\begin{theorem}\label{SpectrumOpMulti}\quad\\
Assume that all wave numbers are equal i.e. $\kappa_{0} = \kappa_{1} = \dots = \kappa_{n}$. 
Let $\mathfrak{S}_{p}(\mrA-\alpha\Pi)$ refer to the point spectrum of $\mrA-\alpha\Pi$
i.e. the set of its eigenvalues. Then the spectrum of this operators coincides with the
point spectrum, and it is given by 
$$
\mathfrak{S}( \mrA-\alpha\Pi ) = \mathfrak{S}_{p}( \mrA-\alpha\Pi ) = \{+\sqrt{1+\alpha^{2}}, -\sqrt{1+\alpha^{2}}\}.
$$
\end{theorem}
\noindent\textbf{Proof:}

The result is clear if $\alpha = 0$ since $(\mrm{Id}+\mrA)/2$ is a projector, so for the remaining of 
this proof, we will assume that $\alpha\neq 0$. Taking the square of the above operator, and 
using Lemma \ref{IdentityDiagOp}, yields $(\mrA-\alpha\Pi)^{2} = (1+\alpha^{2})\mrm{Id} -\alpha 
(\mrm{T}\Pi + \Pi\mrm{T})$. Then it is a direct consequence of Corollary \ref{Nilpotent1} and 
Proposition \ref{Nilpotent2} that the operator $\mrm{T}\Pi + \Pi\mrm{T}$ is nilpotent. Hence
according to \cite[Thm.10.13]{MR1157815} and the spectral mapping theorem \cite[Thm.10.28]{MR1157815}, 
we have $\mathfrak{S}(\mrm{T}\Pi + \Pi\mrm{T}) = \{0\}$. This also shows that
$\mathfrak{S}( (\mrA-\alpha\Pi)^{2}) = \{1+\alpha^{2}\}$, hence applying once again the
spectral mapping theorem, we finally conclude that $\mathfrak{S}(\mrA-\alpha\Pi)\subset 
\{+\sqrt{1+\alpha^{2}}, -\sqrt{1+\alpha^{2}}\}$. 

\quad\\
Denote for a moment $f(\lambda):=\lambda^{2}$. Then clearly $\mathfrak{S}_{p}(\mrm{T}\Pi +
\Pi\mrm{T}) = \{0\}$ since $\mrm{T}\Pi + \Pi\mrm{T}$ is nilpotent. Moreover,
according to the "point spectrum counterpart" of the spectral mapping theorem
\cite[Thm.10.33]{MR1157815}, we have $f(\mathfrak{S}_{p}(\mrA-\alpha\Pi)) =
\mathfrak{S}_{p}(f(\mrA-\alpha\Pi)) = \{1+\alpha^{2}\}$. Hence we conclude
also that $\mathfrak{S}_{p}(\mrA-\alpha\Pi)\subset \{+\sqrt{1+\alpha^{2}}, -\sqrt{1+\alpha^{2}}\}$.
If we can prove that the previous inclusion is actually an equality, then the proof will
be finished. It suffices to show that, if $\lambda\in \mathfrak{S}_{p}(\mrA-\alpha\Pi)$,
then we also have  $-\lambda\in \mathfrak{S}_{p}(\mrA-\alpha\Pi)$.

\quad\\
Take an eigenvector $\ctr{u}\in \mbH(\Sigma)\setminus\{0\}$ of $\mrA-\alpha\Pi$ 
associated to the eigenvalue $\lambda$. Since $\mrA^{2} = \Pi^{2} = \mrm{Id}$, 
we have $\mrA (\mrA\Pi - \Pi\mrA) = -(\mrA\Pi - \Pi\mrA)\mrA$ and 
$\Pi (\mrA\Pi - \Pi\mrA) = -(\mrA\Pi - \Pi\mrA)\Pi$. As a consequence 
$(\mrA-\alpha\Pi) (\mrA\Pi - \Pi\mrA)\ctr{u} = -(\mrA\Pi - \Pi\mrA)(\mrA-\alpha\Pi)\ctr{u}
 = -\lambda (\mrA\Pi - \Pi\mrA)\ctr{u}$. Hence, $-\lambda$ is an eigenvalue of $\mrA-\alpha\Pi$ 
if $\lambda$ is, provided that $(\mrA\Pi - \Pi\mrA)\ctr{u}\neq 0$. Observe that, 
since $(\mrA - \alpha\Pi)\ctr{u} = \lambda\ctr{u}$ and $\mrA^{2} = \Pi^{2} = \mrm{Id}$, 
we have $\Pi\mrA\ctr{u} = \alpha \ctr{u} +\lambda\Pi\ctr{u}$ and $\mrA\Pi\ctr{u} =
(1/\alpha)\ctr{u}-(\lambda/\alpha)\mrA\ctr{u}$. Summing these two identities yields
$$
(\Pi\mrA - \mrA\Pi)\ctr{u} = (\alpha-1/\alpha)\ctr{u}+\frac{\lambda}{\alpha}(\mrA+\alpha\Pi)\ctr{u}.
$$
If we can prove that $0$ is not an eigenvalue of $(\alpha-1/\alpha)\mrm{Id}+\frac{\lambda}{\alpha}(\mrA+\alpha\Pi)$, 
this will show that $(\Pi\mrA - \mrA\Pi)\ctr{u}\neq 0$. From the first part of the proof, we know that the spectrum 
of $\mrA+\alpha\Pi$ is included in $\{+\sqrt{1+\alpha^{2}},-\sqrt{1+\alpha^{2}}\}$. Besides $\lambda$ equals 
$+\sqrt{1+\alpha^{2}}$ or $-\sqrt{1+\alpha^{2}}$. As a consequence the spectrum of 
$(\alpha-1/\alpha)\mrm{Id}+\frac{\lambda}{\alpha}(\mrA+\alpha\Pi)$ only contains the values 
$$
\alpha-\frac{1}{\alpha}\pm\frac{\lambda}{\alpha}\sqrt{1+\alpha^{2}} = 
\frac{\alpha^{2}-1 \pm(1+\alpha^{2})}{\alpha} =\;\; 2\alpha
\;\;\textrm{or}\;\;-2/\alpha
$$
Since $\alpha\in\C\setminus\{0\}$, we have $2\alpha\neq 0$ and $2/\alpha\neq 0$. As a consequence 
the spectrum of $(\alpha-1/\alpha)\mrm{Id}+\frac{\lambda}{\alpha}(\mrA+\alpha\Pi)$ does not contain $0$
so $(\Pi\mrA - \mrA\Pi)\ctr{u}\neq 0$ necessarily. \hfill $\Box$

\quad\\
The previous theorem can be reformulated as follows. 
\begin{corollary}\label{Hyperbole}\quad\\
Assume that all wave numbers are equal $\kappa_{0} = \kappa_{1} = \dots =\kappa_{n}$. Then for any 
pair of complex numbers $\alpha,\beta\in\C$ the operator $\mrA +\alpha\Pi +
\beta\mrm{Id}$ is invertible if and only if $\beta^{2}-\alpha^{2} \neq 1$.
\end{corollary}

\noindent 
Theorem \ref{SpectrumOpMulti} also leads directly to explicit expression for the spectrum of the 
local multi-trace operator. Thanks to Fredholm theory, this implies a well-posedness result.

\begin{corollary}\label{WellPosednessMTFLoc}\quad\\
  For any $\alpha\in \C\setminus\{0\}$, the operator $\mrm{L}_{\alpha} :=
  (1-\alpha)(\mrA-\mrm{Id}) + \alpha(\mrA-\Pi)$ is invertible. Moreover,
  in the case where all wave numbers are equal $\kappa_{0} = \kappa_{1} 
  = \dots = \kappa_{n}$, its spectrum equals its point spectrum and
  $\mathfrak{S}(\mrm{L}_{\alpha}) = \mathfrak{S}_{p}(\mrm{L}_{\alpha}) =
  \{-1+\alpha-\sqrt{1+\alpha^{2}}, -1+\alpha+\sqrt{1+\alpha^{2}}  \}$.
\end{corollary}
\noindent \textbf{Proof:} 

Assume that the wave numbers $\kappa_{0}, \kappa_{1},\dots,\kappa_{n}$ are arbitrary elements 
of $(0,+\infty)$, and consider any $\alpha\in \C\setminus\{0\}$. Let $\mrA_{\star}$ refer to 
the operator defined in the same manner as $\mrA$ but with wave numbers all equal to 
$\kappa_{\star} = \imath$. Then the operator $\mrA - \mrA_{\star}:\mbH(\Sigma)\to\mbH(\Sigma)$ 
is compact as it only involves integral operators with regular kernels, see  
\cite[Lemma 3.9.8]{SauterSchwab}. Then Theorem \ref{SpectrumOpMulti} shows that 
$(1-\alpha)(\mrA_{\star}-\mrm{Id}) + \alpha(\mrA_{\star}-\Pi)$ is invertible, since its 
eigenvalues $-1+\alpha\pm\sqrt{1+\alpha^{2}}$ differ from $0$ as $\alpha\neq 0$. 
Hence $(1-\alpha)(\mrA-\mrm{Id}) + \alpha(\mrA-\Pi)$ is a compact perturbation of an 
isomorphism. According to Fredholm-Riesz-Schauder theory (see \cite[Chap.2]{McLean} for example), 
this operator is invertible if and only if it is one-to-one. Since it is injective according to 
Proposition \ref{Uniqueness}, we finally conclude that $(1-\alpha)(\mrA-\mrm{Id}) + 
\alpha(\mrA-\Pi)$ is an isomorphism. The second statement above concerning the spectrum 
is a trivial consequence of Theorem \ref{SpectrumOpMulti}.
\hfill $\Box$

\quad\\
In the case where wave numbers take arbitrary  values the spectrum is not reduced
to $-1+\alpha\pm\sqrt{1+\alpha^{2}}$ anymore. However a difference of wave numbers
only induces compact perturbation of integral operators so that, in the general
case, this result still indicates the location of accumulation points of the spectrum.

\begin{corollary}\label{WellPosednessMTFLoc2}\quad\\
  For any $\alpha\in \C\setminus\{0\}$, set $\mrm{L}_{\alpha}:= (1-\alpha)(\mrA-\mrm{Id})
  + \alpha(\mrA-\Pi)$. Then any element $\lambda\in \mathfrak{S}(\mrm{L}_{\alpha})\setminus\{
  -1+\alpha+\sqrt{1+\alpha^{2}}, -1+\alpha-\sqrt{1+\alpha^{2}}\}$ is an isolated
  eigenvalue with $\mrm{dim}(\,\mrm{ker}(\mrm{L}_{\alpha} - \lambda\mrm{Id})\,)<+\infty$.
  Moreover the two values $-1+\alpha\pm\sqrt{1+\alpha^{2}}$ are the only possible accumulation
  points of $\mathfrak{S}(\mrm{L}_{\alpha})$.
\end{corollary}
\noindent \textbf{Proof:}

Denote $\mu_{\alpha}^{\pm}:=-1+\alpha\pm\sqrt{1+\alpha^{2}}$, and set $\mathcal{L}(\lambda):=
\mrm{L}_{\alpha} - \lambda\mrm{Id}$. Then $\lambda\mapsto\mathcal{L}(\lambda)$ is an analytic
operator pencil, and it is Fredholm valued for $\lambda\neq \mu_{\alpha}^{\pm}$. Indeed
take any $\lambda\in\C\setminus\{\mu_{\alpha}^{+},\mu_{\alpha}^{-}\}$, and define the
operator $\mrm{L}'_{\alpha}$ in the same manner as $\mrm{L}_{\alpha}$ except
that all wave numbers are taken equal to $\kappa_{0}$. The operator
$\mathcal{L}'(\lambda):= \mrm{L}'_{\alpha} - \lambda\mrm{Id}$ is invertible according to
Corollary \ref{WellPosednessMTFLoc}, and $\mathcal{L}'(\lambda)-\mathcal{L}(\lambda)$ is
compact. As a consequence, $\mathcal{L}(\lambda)$ is a compact perturbation of an invertible
operator, so it is Fredholm with index $0$ and admits finite dimensional kernel.

\quad\\
Further, since $\mrm{L}_{\alpha}$ is a bounded operator, $\mathcal{L}(\lambda)$ is
invertible for $\lambda>\Vert \mrm{L}_{\alpha}\Vert$, where $\Vert\;\;\Vert$ refers
here to the norm naturally associated to continuous operators mapping $\mbH(\Sigma)$
into itself. As a consequence, we can apply Fredholm analytic theorem (see Appendix
A in \cite{MR1729870}) which shows that $\mathcal{L}(\lambda)$ is invertible in $\C
\setminus\{\mu_{\alpha}^{+},\mu_{\alpha}^{-}\}$ except for a countable set of isolated
values. Moreover, we have just seen that all these values can only lie in the
disc of center $0$ and radius $\Vert \mrm{L}_{\alpha}\Vert$. 
\hfill $\Box$

\section{Numerical evidences}\label{Numerics}

In this section, we present a series of numerical results confirming
the conclusions presented previously. We consider 2-D scattering
problems of the form (\ref{InitPb}) involving three domains. As regards discretization,
we consider a uniform paneling $\Sigma^{h}\simeq \Sigma$ which induces a mesh for each
of the sub-domains $\Gamma_{j}^{h}\simeq \Gamma_{j}$, $\Gamma_{j}^{h}\subset \Sigma^{h}$. The
discrete spaces $\mbH_{h}(\Sigma)$ are constructed on these meshes by means of
$\mathbb{P}_{1}$-Lagrange shape functions for both Dirichlet and Neumann traces
\begin{equation}\label{DiscreteSpace}
\begin{array}{l}
\mbH_{h}(\Sigma) = \{\; (u^{h}_{j}, p^{h}_{j})_{j=0,1,2}\;\textrm{such that}\\[5pt]
\hspace{2cm}\forall j = 0,1,2, \textrm{for all panel}\; e\subset \Gamma_{j}^{h},\;\;
u^{h}_{j}\vert_{e},p^{h}_{j}\vert_{e}\in \mathbb{P}_{1}(e)\; \}.
\end{array}
\end{equation}
Denote $\mrm{B}_{h}$ the matrix associated to the Galerkin discretization of the local
multi-trace formulation (\ref{LocMltTr}) by means of the discrete space (\ref{DiscreteSpace}),
and let us denote $\mrm{M}_{h}$ the matrix obtained by Galerkin discretization of the bilinear
form $(\ctru,\ctrv)\mapsto \llbr \ctru,\ctrv\rrbr$. We shall focus our attention
on the spectrum of the matrix $(\mrm{M}_{h})^{-1}\mrm{B}_{h}$ that may be considered as
an approximation of the continuous operator associated to Formulation (\ref{LocMltTr}).

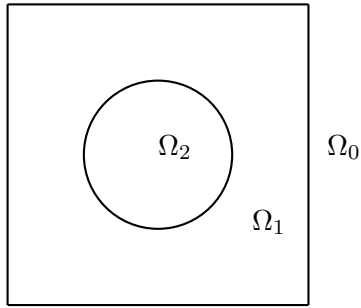
\begin{figure}[h]
\begin{pspicture}(0,4.25)(0,0)

  \put(5.5,-0.75){
    \put(0,1){
      \psline(0,0)(4,0)
      \psline(4,0)(4,4)
      \psline(4,4)(0,4)
      \psline(0,4)(0,0)    
      \pscircle(2,2){1} }
    
    \put(4.25,3){$\Omega_{0}$}
    \put(3.25,2){$\Omega_{1}$}
    \put(2,3){$\Omega_{2}$}
  }
  
\end{pspicture}
\caption{First geometrical configuration} \label{Fig1}
\end{figure}

\noindent 
All computations have been achieved on a laptop equipped with a 2-core Intel i7-3520M processor 
at 2.9GHz with 4 GB of RAM. Meshes have been generated using Gmsh \cite{MR2566786} (see also 
the website \verb?http://geuz.org/gmsh/?). For computation of eigenvalues we relied on the
Arpack++ library  (see \verb?http://www.ime.unicamp.br/~chico/arpack++/?).

\begin{figure}[h]
\begin{pspicture}(0,4.5)(0,0)
  \put(0,0){\includegraphics[width=7cm]{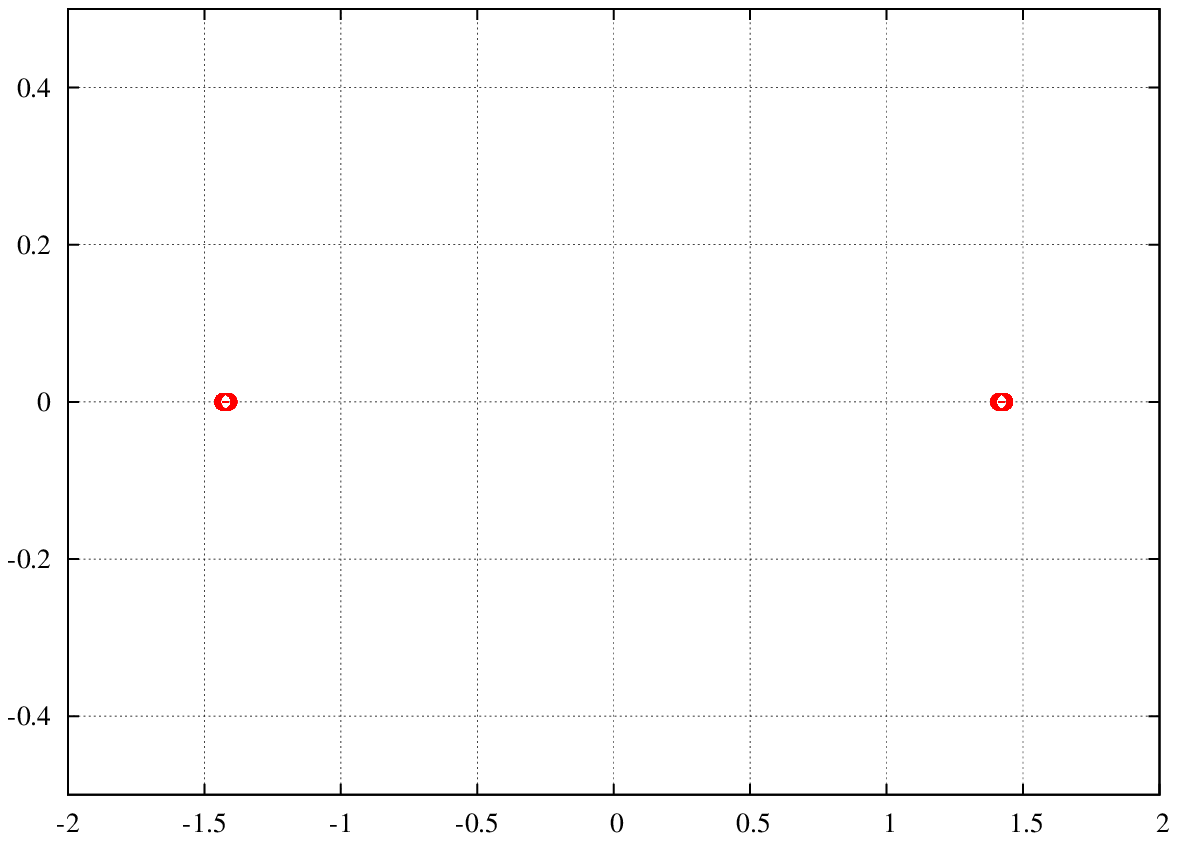}}
  \put(8,0){\includegraphics[width=7cm]{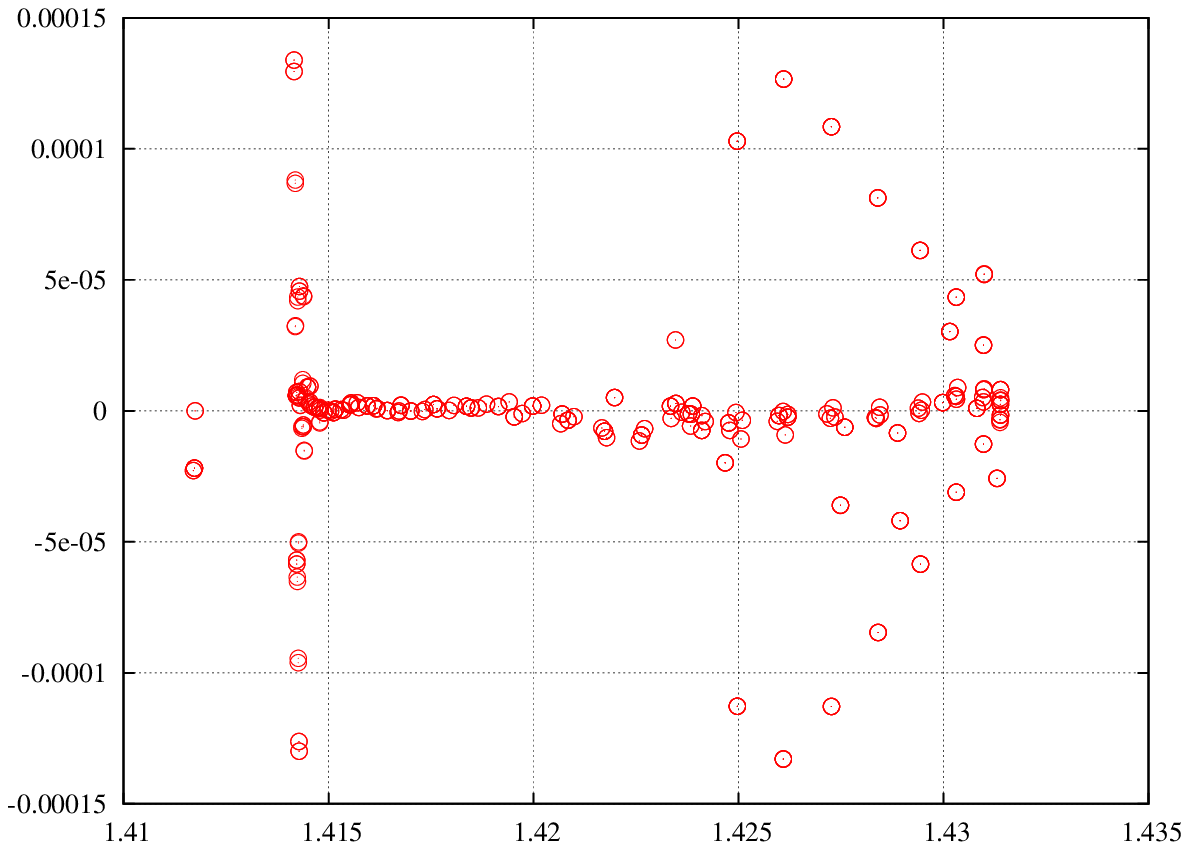}}
  \pspolygon[linestyle=dashed,linewidth=0.5pt,dash=2pt 2pt](5.6,2.3)(6.1,2.3)(6.1,2.7)(5.6,2.7)
  \psline[linestyle=dashed,linewidth=0.5pt,dash=2pt 2pt](6.1,2.7)(8.75,4.75)
  \psline[linestyle=dashed,linewidth=0.5pt,dash=2pt 2pt](6.1,2.3)(8.75,0.25)
  
\end{pspicture}
\caption{Two clusters of eigenvalues in the case of no contrast of wave numbers (left)
  and a zoom on the cluster of positive eigenvalues (right). We took
  $\alpha = 1$ and $\kappa_{0} = \kappa_{1} = \kappa_{2} = 1$.} \label{Fig2} 
\end{figure}

\quad\\[-5pt]
Figure \ref{Fig1}  represents a first geometrical configuration. The boundary of
$\Omega_{0}$ is a unit square centered at $0$, and the boundary of $\Omega_{2}$ is a
circle of radius $0.5$ centered at $0$. Figure \ref{Fig2} represents the spectrum of
$(\mrm{M}_{h})^{-1}\mrm{B}_{h}$ for a mesh width $h = 0.05$ and $\kappa_{0} = \kappa_{1} =
\kappa_{2} = 1$.

\quad\\

\begin{figure}[!h]
\begin{pspicture}(0,4.5)(0,0)
  \put(0,0){\includegraphics[width=7cm]{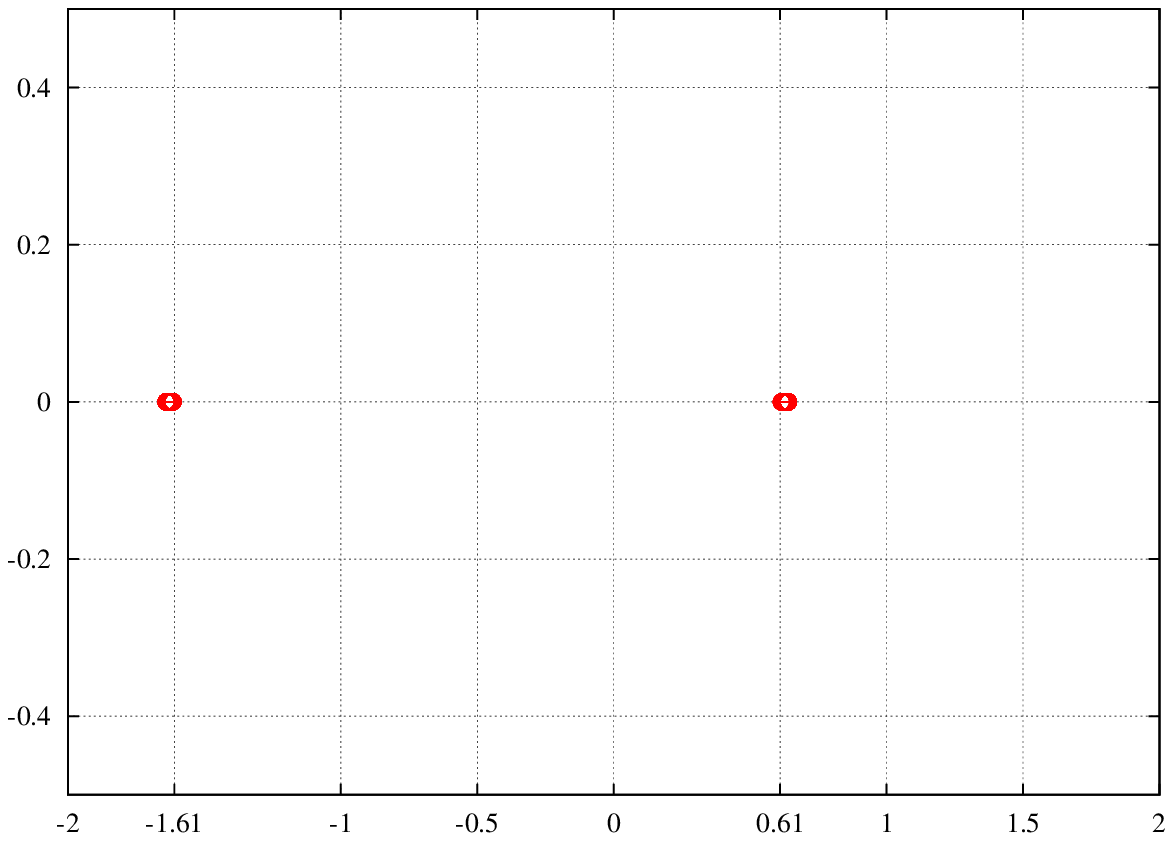}}
  \put(8,0){\includegraphics[width=7cm]{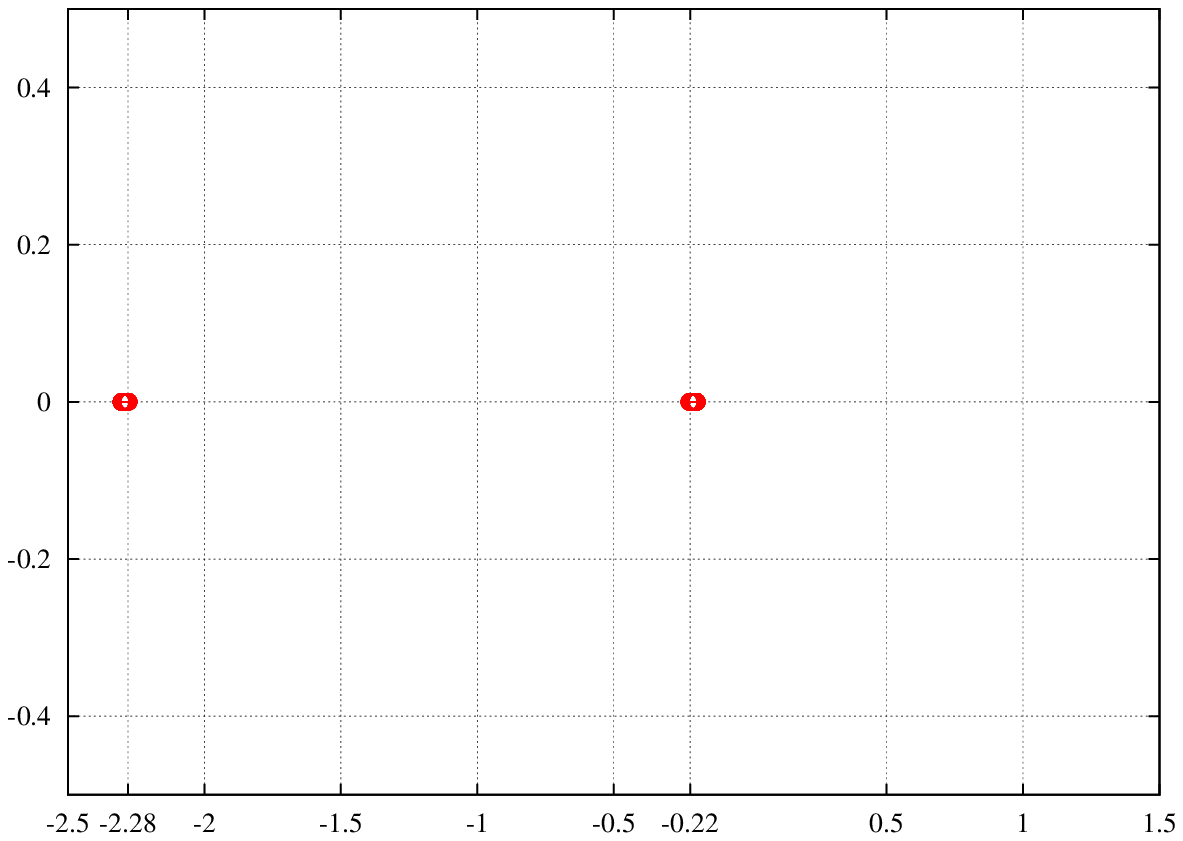}}
  
\end{pspicture}
\caption{Spectrum of the local multi-trace operator
  with $\kappa_{0} = \kappa_{1} = \kappa_{2} = 1$  and
  $\alpha = 0.5$ (left) or $\alpha = -0.25$ (right) } \label{Fig3} 
\end{figure}

\quad\\[-5pt]
The spectrum clearly takes the form of two clusters centered at the values
$\pm\sqrt{2} = -1+\alpha\pm\sqrt{1+\alpha^{2}}$ for $\alpha = 1$. Figure \ref{Fig3} shows
the spectrum of the same matrix, with the same geometrical
configuration, but with $\alpha = 0.5$ and $\alpha = -0.25$. The formula $-1+\alpha\pm\sqrt{1+\alpha^{2}}$
yields the values $0.61803$ and $-1.6180$ for $\alpha = 0.5$ (up to $5$ digits),
and $-0.21922$ and $-2.2808$ for  $\alpha = -0.25$, which is consistent with our theory.

\quad\\

\begin{figure}[h]
\begin{pspicture}(0,4)(0,0)
  \put(4,0){\includegraphics[width=7cm]{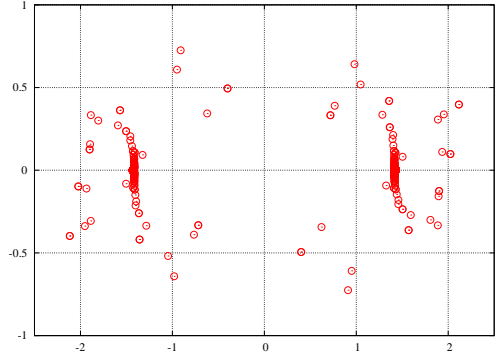}}
\end{pspicture}
\caption{Spectrum of the local multi-trace operator in the case
  $\kappa_{0} = 1$, $\kappa_{1} = 5$ and $\kappa_{2} = 2$ with $\alpha = 1$. } \label{Fig4} 
\end{figure}

\quad\\
Next, in Figure \ref{Fig4}, we consider the case $\alpha = 1$, but wave numbers differ taking the values
$\kappa_{0} = 1$, $\kappa_{1} = 5$ and $\kappa_{2} = 2$. The mesh width remains $h = 0.05$.
Although the eigenvalues are not clustered anymore, they are more densely grouped around $\pm\sqrt{2}$
suggesting that these are the only two accumulation points of the spectrum of the continuous operator.

\quad\\

\noindent
For the next series of figures, we consider the same scattering problem, but
in a different geometrical configuration. The new configuration is depicted in
the picture below: there are two square scatterer separated by a thin gap of
width delta. In the theory we have presented,  we needed to assume that there is no junction
point i.e. points where at least three sub-domains abut. We wish to test the robustness,
with respect to this assumption, of the theoretical formulas obtained.

\begin{pspicture}(0,4.5)(0,0)
    
\put(7,0.5){  
\pspolygon(-3,0)(-0.25,0)(-0.25,3)(-3,3)
\pspolygon(0.25,0)(3,0)(3,3)(0.25,3)
\psline[linestyle=dashed](0.25,0)(0.25,-0.75)
\psline[linestyle=dashed](-0.25,0)(-0.25,-0.75)
\psline[linestyle=dashed,arrowinset=0,arrowsize=4pt]{<->}(-0.25,-0.5)(-0.75,-0.5)(+0.75,-0.5)(+0.25,-0.5)
\put(-0.05,-0.2){$\delta$}

\put(-1.65,1.5){$\Omega_{1}$}
\put(+1.25,1.5){$\Omega_{2}$}
\put(+3.5,1.5) {$\Omega_{0}$}
\pscircle[fillstyle=solid,fillcolor=black](3,0){0.1}
\pscircle[fillstyle=solid,fillcolor=black](3,3){0.1}
\pscircle[fillstyle=solid,fillcolor=black](-3,3){0.1}
\pscircle[fillstyle=solid,fillcolor=black](-3,0){0.1}
\put(2.25,-0.45){\tiny $(1,-0.5)$}
\put(2.25,3.25){\tiny $(1,+0.5)$}
\put(-3.25,3.25){\tiny $(-1,+0.5)$}
\put(-3.25,-0.45){\tiny $(-1,-0.5)$} }

\end{pspicture}

\begin{figure}[!h]
\begin{pspicture}(0,3.25)(0,0)
  \put(0,0.5){\includegraphics[width=5cm]   {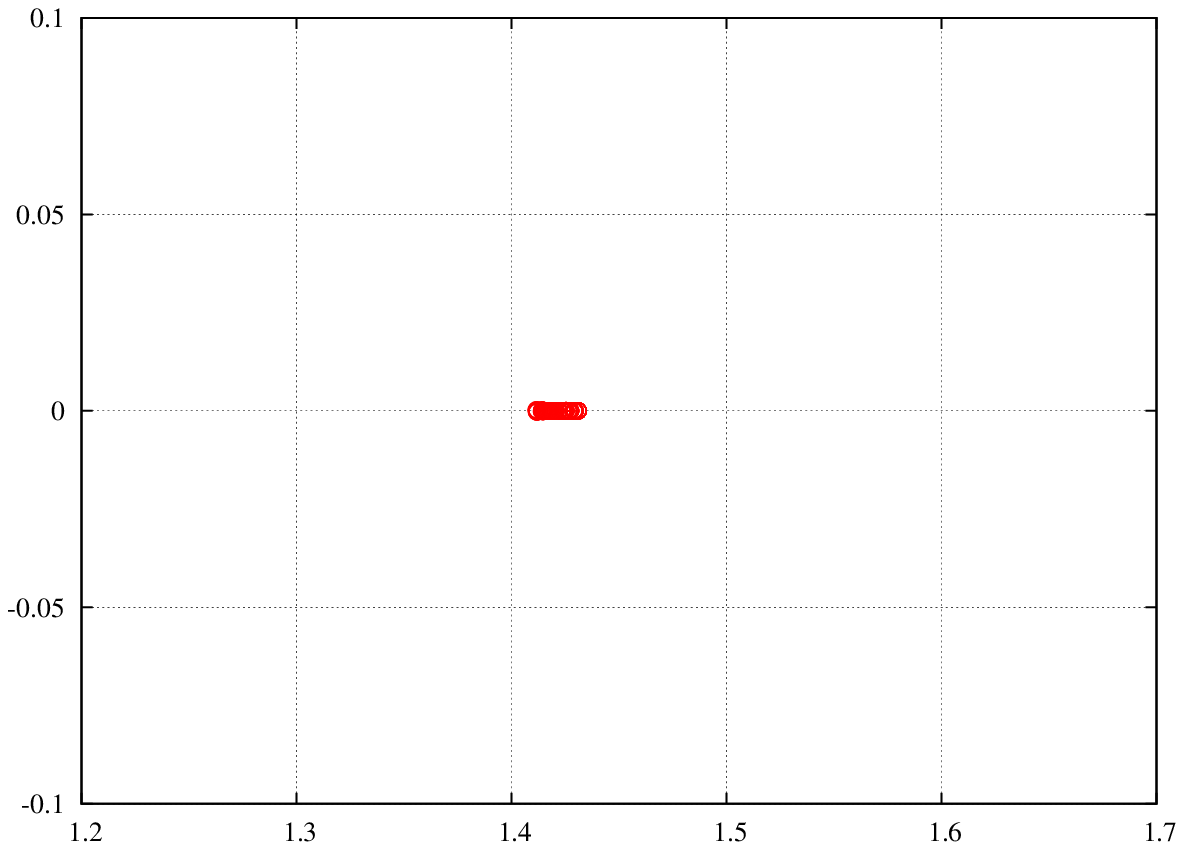}}
  \put(5.25,0.5){\includegraphics[width=5cm]{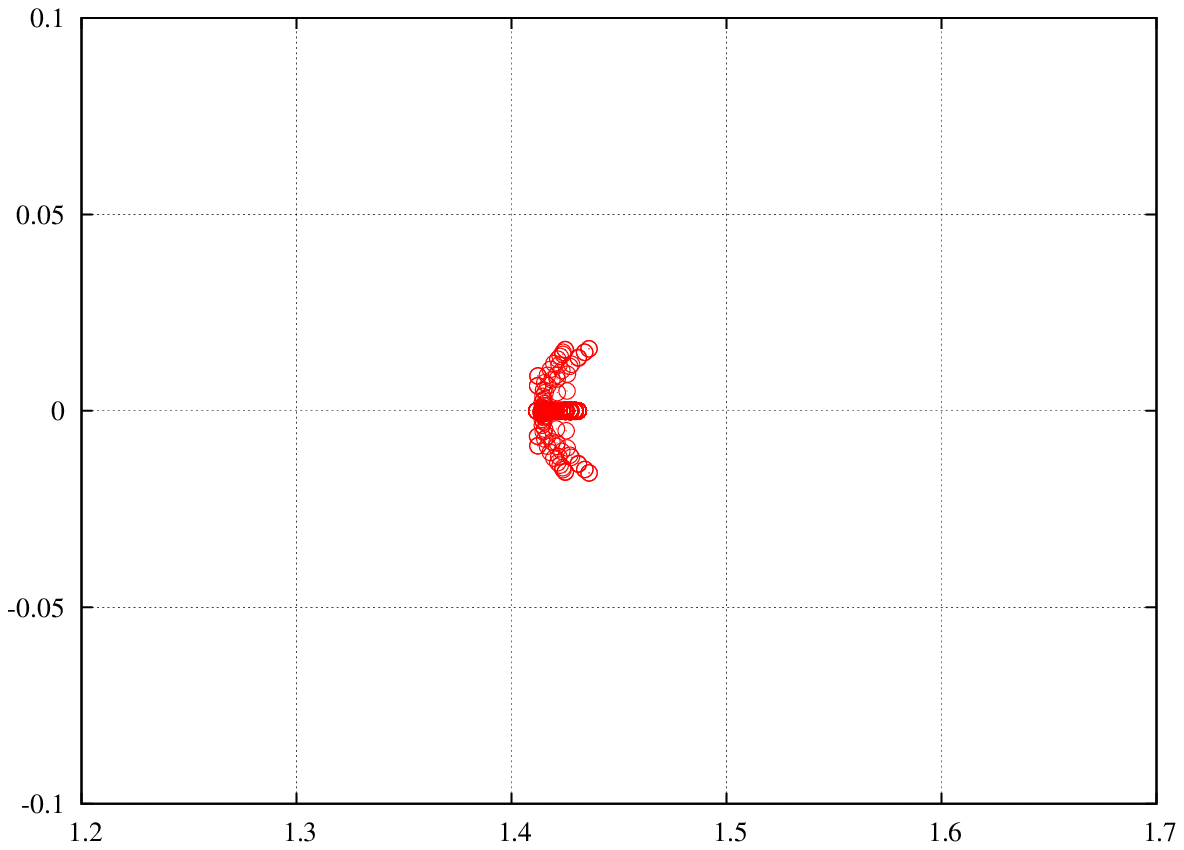}}
  \put(10.5,0.5){\includegraphics[width=5cm]{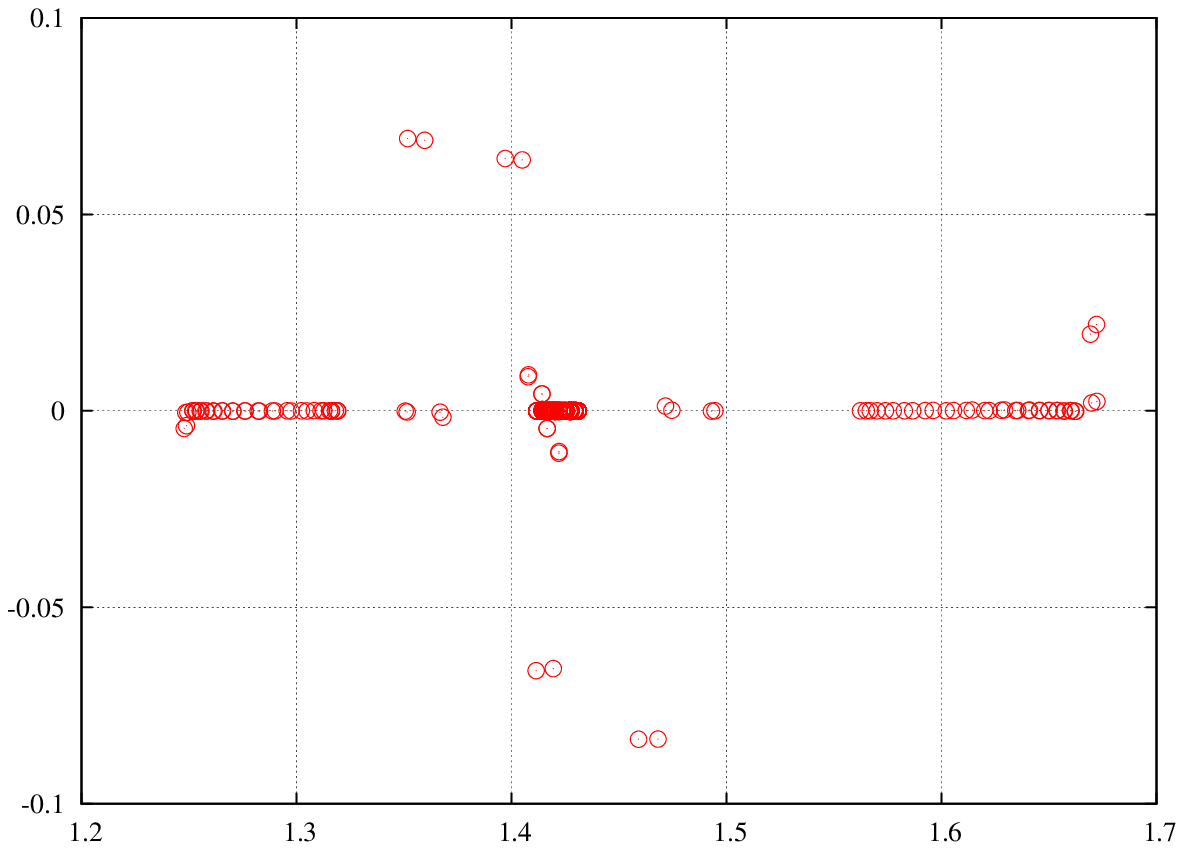}}

  \put(2,0)   {\textbf{(a)}}
  \put(7.25,0){\textbf{(b)}}
  \put(12.5,0){\textbf{(c)}}

\end{pspicture}
\caption{ Spectrum of the local multi-trace operator for
  $\kappa_{0} = \kappa_{1} = \kappa_{2} = 1$, $\alpha = 1$ and three
  different values of gap: $\delta = 0.1$ (left) $\delta = 0.01$ (center) and
  $\delta = 0.001$ (right) } \label{Fig5} 
\end{figure}

\quad\\
In Figure \ref{Fig5}, we consider $\kappa_{0} = \kappa_{1} = \kappa_{2}$ and
$\alpha = 1$. For a fixed strictly positive value of $\delta>0$, the eigenvalues
are clustered around $\pm\sqrt{2}$. Each of the figures below represent a zoom
at the cluster centered at $+\sqrt{2}$ for various values of $\delta$.

\quad\\

\quad\\
The cluster is more and more scattered as
the gap closes, suggesting  that the assumption that there is no junction
point is mandatory, in spite of quadrature rules being less reliable as $\delta$
is close to $0$. In the last picture below, we examine the case where the gap is
closed i.e. $\delta = 0$, which corresponds to the presence of a junction point in
the geometry.

\quad\\
\begin{figure}[h]
\begin{pspicture}(0,4.5)(0,0)
  \put(0.25,0.){\includegraphics[width=7cm]{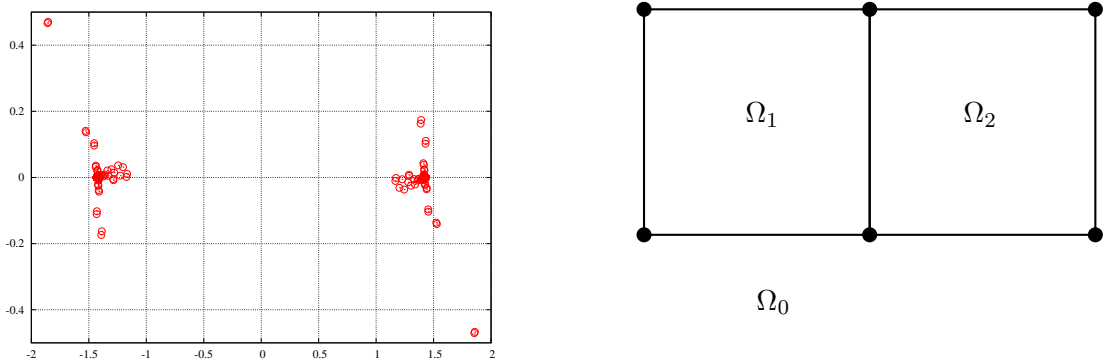}}

  \put(12,1.75){
  \pspolygon(-3,0)(0.,0)(0.,3)(-3,3)
  \pspolygon(0.,0)(3,0)(3,3)(0.,3)
  \pscircle[fillstyle=solid,fillcolor=black](3,0){0.1}
  \pscircle[fillstyle=solid,fillcolor=black](3,3){0.1}
  \pscircle[fillstyle=solid,fillcolor=black](0,0){0.1}
  \pscircle[fillstyle=solid,fillcolor=black](0,3){0.1}
  \pscircle[fillstyle=solid,fillcolor=black](-3,3){0.1}
  \pscircle[fillstyle=solid,fillcolor=black](-3,0){0.1}
  \put(-1.65,1.5){$\Omega_{1}$}
  \put(+1.25,1.5){$\Omega_{2}$}
  \put(-1.5,-1) {$\Omega_{0}$}

  }

\end{pspicture}
\caption{ Spectrum of the local multi-trace operator (left) for
  $\kappa_{0} = \kappa_{1} = \kappa_{2} = 1$, $\alpha = 1$
  in the presence of junction points in the geometry (right).} \label{Fig6} 
\end{figure}

\quad\\[-5pt]
Unfortunately, geometrical configurations involving junction points are not covered
by the theory of the present article. However the result of Figure \ref{Fig6} suggests that,
although the eigenvalues are not anymore closely clustered around the values
$-1+\alpha\pm\sqrt{1+\alpha^{2}}$, these two theoretical points remain accumulation points
of the spectrum of the local multi-trace operator.

\paragraph{Acknowledgement}
The author would like to thank V.Dolean and M.J.Gander for fruitful discussions.

\bibliographystyle{plain}
\def\cprime{$'$}

\end{document}